%% file: brief.tex
\documentclass[journal]{IEEEtran}
\usepackage{cite}
\usepackage{amsmath,amssymb,amsfonts, amsthm}
\usepackage{algorithm}
\usepackage{algpseudocode}
\usepackage{graphicx}
\usepackage{textcomp}
\usepackage{xcolor}
\usepackage{acro}
\usepackage{siunitx}
\usepackage{tikz,pgfplots}
\usepackage{hyperref}
\usepackage{setspace}
\usepackage{caption}
\usepackage{subcaption}
\allowdisplaybreaks

\usetikzlibrary{positioning}
\usetikzlibrary{calc}
\usetikzlibrary{arrows}
\usetikzlibrary{shapes.geometric}

\newtheorem{proposition}{Proposition}

\def\BibTeX{{\rm B\kern-.05em{\sc i\kern-.025em b}\kern-.08em
		T\kern-.1667em\lower.7ex\hbox{E}\kern-.125emX}}
\newcommand{\dc}{\mathrm{dc}}
\newcommand{\abc}{\mathrm{abc}}
\newcommand{\A}{\mathrm{a}}
\newcommand{\B}{\mathrm{b}}
\newcommand{\C}{\mathrm{c}}
\newcommand{\g}{\mathrm{g}}
\newcommand{\ph}{\mathrm{ }} % dropped ph index because ambiguous
\newcommand{\sw}{\mathrm{sw}}
\newcommand{\s}{\mathrm{s}}

\newcommand{\p}{\mathrm{p}}
\newcommand{\U}{\mathrm{u}}

\newcommand{\Sim}{\mathrm{sim}}
\newcommand{\I}{\b{I}}
\newcommand{\0}{\b{0}}

\newcommand{\ini}{\mathrm{ini}}
\newcommand{\augp}{\mathrm{T}}
\newcommand{\augs}{\mathrm{S}}
\newcommand{\opt}{\mathrm{opt}}
\newcommand{\diag}{\mathrm{diag}}
\newcommand{\prev}{\mathrm{prev}}

\newcommand{\ct}{\mathrm{FT}}
\newcommand{\cs}{\mathrm{FL}}
\newcommand{\past}{\mathrm{past}}

\definecolor{blueCol}{rgb}{0, 0.447 0.741}

\renewcommand{\b}[1]{\boldsymbol{#1}}

\renewcommand{\th}{\mathrm{th}}

\newlength{\fheight}
\newlength{\fwidth}

\DeclareSIUnit{\pu}{pu}

\DeclareAcronym{FCS}{short = FCS, long = finite control set}
\DeclareAcronym{MPC}{short = MPC, long = model predictive control}
\DeclareAcronym{PI}{short = PI, long = proportional-integral}
\DeclareAcronym{MIMO}{short = MIMO, long = multiple-input multiple-output}
\DeclareAcronym{ILS}{short = ILS, long = integer least squares}
\DeclareAcronym{3L-NPC}{short = 3L-NPC, long = 3-level neutral-point-clampled}
\DeclareAcronym{IIR}{short = IIR, long = infinite impulse response}
\DeclareAcronym{TDD}{short = TDD, long = total demand distortion}
\DeclareAcronym{PLL}{short = PLL, long = phase-locked-loop}
\DeclareAcronym{IQP}{short = IQP, long = integer-quadratic program}
\DeclareAcronym{c1}{short = PI-MPC, long = hierarchical controller with a conventional FCS-MPC on the inner loop and a proportional-integral controller for switching frequency tracking on an outer loop}
\DeclareAcronym{c2}{short = FT-MPC, long = switching frequency tracking FCS-MPC}
\DeclareAcronym{c3}{short = FL-MPC, long = switching frequency limiting FCS-MPC}
\DeclareAcronym{sf}{short = switching frequency, long = switching frequency}

\begin{document}
	\title{Switching Frequency Limitation with Finite Control Set Model Predictive Control via Slack Variables}
	\author{Luca M. Hartmann, Orcun Karaca, Tinus Dorfling, Tobias Geyer
		\thanks{ Corresponding author: Orcun Karaca.
L. M. Hartmann is with ETH Z{\"u}rich, Switzerland. email: {\tt hartmannl@ethz.ch.} O. Karaca and T. Dorfling are with ABB Corporate Research Center, Switzerland. emails: {\tt \{orcun.karaca,} {\tt martinus-david.dorfling\}@abb.ch.com}. T. Geyer is with ABB System Drives, Switzerland. email: {\tt t.geyer@ieee.org.} This work was part of the thesis work of L. M. Hartmann at ETH Z{\"u}rich. The authors thank Prof. F. D{\"o}rfler of Automatic Control Laboratory for his feedback and this opportunity.}}
	
	\maketitle
	
	\begin{abstract}
	Past work proposed an extension to finite control set model predictive control to simultaneously track both a current reference and a switching frequency reference. Such an objective can jeopardize the current tracking performance, and this can potentially be alleviated by instead limiting the switching frequency. To this end, we propose to limit the switching frequency in finite control set model predictive control. The switching frequency is captured with an infinite impulse response filter and bounded by an inequality constraint; its corresponding slack variable is penalized in the cost function. To solve the problem efficiently, a sphere decoder with a computational speed-up is presented. \looseness=-1
	\end{abstract}
	

\vspace{-.1cm}
\section{Introduction}
\label{sec:introduction}
Given the developments in mathematical optimization techniques and their applications on embedded systems, model predictive control (MPC) has established itself as a promising control methodology in power conversion systems, where the system time constants are well below $1$ms~\cite{karamanakos2020model,rodriguez2021latest,rodriguez2021latest2,harbi2023model,zafra2023long}. Several different variants of MPC have been developed for converter control. Among those, finite control set MPC (FCS-MPC) has gained a lot of attention for its advantages, such as its intuitive design, simple implementation, and high dynamic performance. \looseness=-1

FCS-MPC achieves regulation of the states along their references by directly manipulating the switch positions of a power electronic converter. Since the switch positions are discrete, the resulting optimization problem is an integer program. Its first variants utilized a horizon of one step and an exhaustive enumeration to solve the underlying integer program~\cite{muller2003new,rodriguez2004predictive}. Whenever a linear prediction model is available, an efficient branch-and-bound algorithm called the sphere decoder can be adopted~\cite{geyer_2014, geyer_2016, dorfling_2020}. This enables the use of long horizons within short sampling intervals, which can bring performance benefits to different MPC variants, e.g., increasing the closed-loop stability margin~\cite{karamanakos_2014,geyer2014performance,bordons2015basic}. \looseness=-1

One of the main shortcomings of FCS-MPC is the tuning of its variable switching frequency, which is strongly correlated with the switching losses of the semiconductors~\cite{Rajapakse_2005, Al-Naseem_2000}. Most commonly, the control effort is penalized with the aim to reduce the average switching frequency. 
The tuning of this penalty to attain a particular average switching frequency is known to be sensitive to the other control parameters and the operation point, which greatly complicates the tuning task. Without any switching penalty, on the other hand, the switching frequency would be limited only by the choice of the sampling interval~\cite{karamanakos2019guidelines} (see \cite[Ch. 4]{rodriguez2012predictive} for an empirical study), rendering it similar to a deadbeat-type controller. 
To avoid the tuning process altogether, alternative MPC formulations have been proposed to set the switching frequency directly by fixing the number of possible switching transitions in a given time interval, see~\cite{tarisciotti2014modulated,vazquez2014predictive,karamanakos2018fixed,yang_2021,karamanakos_2021} and the references therein. However, these modified controllers differ significantly from the original FCS-MPC method, and thus they are considered out of scope for this paper. \looseness=-1

To control the switching frequency of FCS-MPC, the authors in~\cite{stellato_2017}~have shown that the penalty term can be replaced with a more meaningful term: a switching frequency tracking term. A second-order infinite impulse response filter can be utilized to capture the predicted switching frequency throughout the horizon. Even more importantly, this modification does not imply any major changes to the FCS-MPC formulation. \looseness=-1

The current tracking performance generally improves with an increase in the switching frequency. Yet, this general trend does not necessarily hold for all cases, for example, there exist certain time instant, in particular during transients, where a momentary decrease in switching frequency could help improve current tracking.
In such cases, when using the formulation proposed in~\cite{stellato_2017}, the two tracking objectives contradict each other. However, a switching frequency lower than a certain limit is in fact not a concern, since it is used as a surrogate to limit the switching losses or even the device temperatures. To this end, the main contribution of this paper is to formulate an FCS-MPC method that limits the switching frequency. This will be achieved by penalizing the slack variable corresponding to a constraint {that upper bounds} the switching frequency. \looseness=-1

Solving FCS-MPC with hard or soft state constraints is in fact rather straightforward in the case of exhaustive enumeration. However, this is largely an unsolved problem when using the sphere decoder. As mentioned before, the sphere decoder is only applicable to optimization problems with linear prediction models. The prediction model for the slack variable of the switching frequency constraint has a nonlinear input-output relation. This nonlinearity is due to a max operator in its output function and an absolute value operator in its state update. One of the first attempts to include constraints (in particular, current constraints) in the sphere decoder has been made in \cite{karamanakos2016constrained,rossi2022constrained}. The authors have proposed an algorithm to map the state constraints to the input constraints, if the state itself has a linear prediction model. A more recent work by~\cite{keusch2023long} took an entirely different perspective by working with linear Gaussian models and constraints expressed in terms of Gaussion priors. However, this approach comes with a loss in optimality. Instead, the approach in this paper will augment the input vector with the slack variables to then formulate the corresponding sphere {decoding} algorithm, resembling the approach in~\cite{liegmann_2017} that included the neutral-point dynamics in FCS-MPC. Moreover, we will provide a computational speed-up to this sphere decoder by exploiting a provable lower bound that incorporates the prediction model of the slack variable we are interested in. \looseness=-1

Our contributions are as follows. We propose an FCS-MPC method that limits the switching frequency by penalizing slack variable of a constraint on the switching frequency.
We formulate a sphere decoder that utilizes slack variables. Finally, a significant computational speed-up is achieved by exploiting a lower bound on the future cost to be incurred from the slack variables corresponding to the switching frequency constraints. \looseness=-1

\vspace{-.1cm}
\section{Preliminaries}\label{sec:system_modeling}
\vspace{0cm}
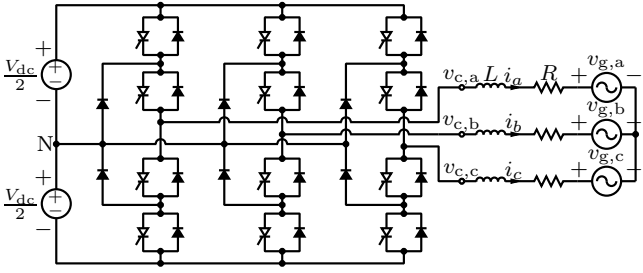
\begin{figure}[t]
    \centering
     \resizebox{.985\columnwidth}{!}{\input{Figures/npc.tex}}
    \caption{Grid-connected 3L-NPC converter.}
    \label{fig:3L-NPC_scheme}
\end{figure}
Consider a \ac{3L-NPC} converter connected to the grid, as depicted in Figure \ref{fig:3L-NPC_scheme}.\footnote{The control methods under consideration are applicable to any system governed by linearized dynamics (e.g., induction machines as in \cite{stellato_2017}).} The half dc-link voltages $\frac{V_\dc}{2}$ are realized by ideal voltage sources, and the neutral-point-potential is fixed to zero; its balancing is out of scope, we kindly refer to \cite{liegmann_2017, Jin_2022}. 
The grid is modeled with an RL-load in series with an infinite bus voltage. Assume all grid parameters and states are known. %Other than having a grid estimation algorithm, this assumption also holds when there is a strong grid connection with the RL-load representing a known transformer impedance.  

Here, the converter voltages are denoted by $\b{v}_{\C} = [v_{\C,\A}\ v_{\C,\B}\ v_{\C,\C}]^\top \in \{-\frac{V_\dc}{2},\ 0,\ \frac{V_\dc}{2}\}^3$. The vectors $\b{i}_\abc = [i_\A\ i_\B \ i_\C]^\top$ and $\b{v}_\g = [v_{\g,\A}\ v_{\g,\B}\ v_{\g,\C}]^\top$ represent the current and the grid voltage, respectively. Three-phase variables in the abc-coordinates are transformed into $\b{\xi}_{\alpha\beta}$ in the stationary $\alpha\beta$ reference frame by $\b{\xi}_{\alpha\beta} = \b{K\xi}_{\A\B\C}$. The inverse operation is denoted by $\b{\xi}_{\A\B\C} = \b{K}^{-1}\b{\xi}_{\alpha\beta}$. The matrices $\b{K}$ and $\b{K}^{-1}$ are defined by the Clarke transformation and its pseudo-inverse, $\b{K} = \frac{2}{3}{\tiny\begin{bmatrix}
			1 & -\frac{1}{2} & -\frac{1}{2} \\
			0 & \frac{\sqrt{3}}{2} & -\frac{\sqrt{3}}{2}
		\end{bmatrix}}$, and $ 
		\b{K}^{-1} = \frac{3}{2}\b{K}^\top.$
	The subscript $\alpha\beta$ is dropped from $\b{\xi}_{\alpha\beta}(k)$, unless stated otherwise.

\vspace{0cm}
\subsection{Physical system model}
\vspace{0cm}
The continuous-time current dynamics for Figure \ref{fig:3L-NPC_scheme} are: $$\frac{d \b{i}_{\alpha\beta}(t)}{dt} = -\frac{R}{L}\b{I}_2\b{i}_{\alpha\beta}(t) + \frac{1}{L}\b{K}\b{v}_{\C,\abc}(t) - \frac{1}{L}\b{K}\b{v}_\g(t),$$ 
with the 2$\times$2 identity matrix $\b{I}_2$, the converter voltage $\b{v}_{\C,\abc}(t) = \frac{V_\dc}{2}\b{u}_{\ph}(t)$, and the switch position $\b{u}_\ph(t)=[{u}_a(t)\ {u}_b(t)\ {u}_c(t)]^\top\in \{-1,0,1\}^3$. Using exact Euler discretization, the discrete-time model becomes
\begin{equation}\label{eq:0}
\begin{split}
    &\b{x}_\ph(k+1) = \b{A}_\ph\b{x}_\ph(k) + \b{B}_\ph\b{u}_\ph(k) + \b{D}_\ph\b{v}_\g(k), \\
    &\b{y}_\ph(k) = \b{C}_\ph \b{x}_\ph(k),
    \end{split}
\end{equation}
where $\b{x}_\ph = \b{i}_{\alpha\beta}$ is the state, $\b{y}_\ph$ the output, $T_\s$ the control sampling interval, and $k\in \mathbb{N}$ the discrete-time index. \looseness=-1
	
\subsection{Switching frequency estimation}
The switching frequency at time step $k$, $f_\sw(k)$, is defined as the number of on- or off-transitions averaged over the  number of semiconductor devices and divided by the sampling interval. As discussed in \cite{stellato_2017}, a discrete-time second-order \ac{IIR} filter can be utilized to estimate the device switching frequency: \looseness=-1
\begin{equation}\label{eq:5}
\begin{split}
    &\resizebox{.91\columnwidth}{!}{%
    $\b{x}_\sw(k+1)
    = \underbrace{ \begin{bmatrix}
            a_1 & 1-a_1 \\
            0 & a_2
        \end{bmatrix}^\top }_{\mbox{\large $\b{A}_\sw$}} \b{x}_\sw(k) + \underbrace{ \frac{1-a_2}{12T_s}\begin{bmatrix}
            1 & 0 \\
            1 & 0 \\
            1 & 0
       \end{bmatrix}^\top }_{\mbox{\large $\b{B}_\sw$}} \b{p}(k),$%
    } \\
    &y_\sw(k) 
    = f_\sw(k) = [0\ 1]\b{x}_\sw(k) = \b{C}_\sw\b{x}_\sw(k). 
     \end{split}
\end{equation}
Here, $a_1, a_2 \in [0, 1]$ are the tuning parameters, and $\b{p}(k) = |\Delta \b{u}_\ph(k)|$ is the element-wise absolute value of the three-phase switching transition $\Delta \b{u}_\ph(k) = \b{u}_\ph(k) - \b{u}_\ph(k-1)$. The filtering window increases and the bandwidth reduces as $a_1$ and $a_2$ are set closer to 1, see~\cite{oppenheim1997signals}. \looseness=-1

\subsection{FCS-MPC for $f_\sw$ tracking}
\begin{figure}[t!]    
    \centering
    \resizebox{.985\columnwidth}{!}{
    \input{Figures/block_diagram.tex}
    }
    \caption{Block diagram of FCS-MPC with switching frequency estimation (EST).}
    \label{fig:MPC-tracking_block-diag}
\end{figure}
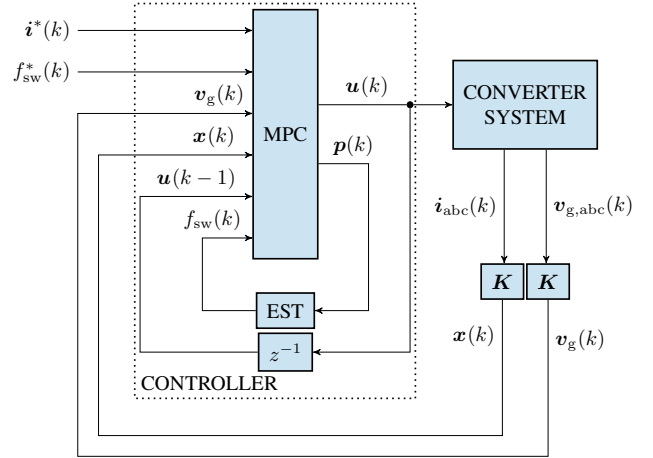

The \ac{c2} of~\cite{stellato_2017} is an MPC algorithm with a prediction horizon of length $N_\p \in \mathbb{N}$, following the block diagram representation in Figure~\ref{fig:MPC-tracking_block-diag}. The \ac{c2} differs from the conventional FCS-MPC in~\cite{geyer_2016} by tracking both the current reference and the~switching frequency reference. 

Define an augmented system with state  $\b{x}_\augp(k) = [\b{x}_\ph(k)^\top\ \b{x}_\sw(k)^\top]^\top$, input $\b{u}_\augp(k) = [\b{u}_\ph(k)^\top\ \b{p}(k)^\top]^\top$, and output $\b{y}_\augp(k) = [\b{y}_\ph(k)^\top\ y_\sw(k)]^\top$. The resulting dynamics are $\b{x}_\augp(k+1) = \b{A}_\augp\b{x}_\augp(k) + \b{B}_\augp\b{u}_\augp(k) + \b{D}_\augp\b{v}_g(k)$, $\b{y}_\augp(k+1) = \b{C}_\augp \b{x}_\augp(k+1)$, where $\b{A}_\augp$, $\b{B}_\augp$, $\b{C}_\augp$, and $\b{D}_\augp$ are defined according to~\eqref{eq:0} and~\eqref{eq:5}. In the remainder, for any vector $\b{\xi}(k)$, we also let $\b{\Xi}(k) = [\b{\xi}(k)^\top\ \b{\xi}(k+1)^\top\ \ldots\ \b{\xi}(k+N_\p - 1)^\top ]^\top$ denote its full-horizon vector for the horizon $N_\p$.\looseness=-1

\ac{c2} solves the following optimization problem at each time step $k$:\looseness=-1
{\medmuskip=2.2mu\begin{equation}\label{eq:c2-opt}
   \begin{split}
    &\min \hspace{-.15cm} \sum_{\ell=k}^{k+N_\p-1} || \b{y}_\augp^*(\ell+1) - \b{y}_\augp(\ell+1) ||_{\b{Q}_\augp}^2 + \lambda_\U || \Delta \b{u}_\ph(\ell) ||_2^2  \\
    &\
        \mathrm{s.t.}\ \b{x}_\augp(\ell+1) = \b{A}_\augp\b{x}_\augp(\ell) + \b{B}_\augp\b{u}_\augp(\ell) + \b{D}_\augp\b{v}_g(\ell),\forall \ell, \\
     & \quad\quad   \b{y}_\augp(\ell+1) = \b{C}_\augp \b{x}_\augp(\ell+1), \forall \ell, \\
     & \quad\quad   \b{u}_\ph(\ell) \in \{-1,\ 0,\ 1\}^{3}, \forall \ell, 
    \end{split}
\end{equation}}where $||\b{\Xi}||_{\b{Q}}^2 = \b{\Xi}^\top\b{Q\Xi}$ refers to the Mahalanobis norm for a given $\b{Q}\succ 0$, $\b{\xi}^*$ denotes the reference of a vector $\b{\xi}${, and $\ell\in\mathbb{N}$ is used to iterate over the prediction horizon}. We have $\b{Q}_\augp = \diag([1\ 1\ \lambda_\sw])$, which uses the relative~switching frequency penalty $\lambda_\sw > 0$ as a trade off between current and switching frequency tracking. The switching penalty $\lambda_\U > 0$ is set to be very small to not interfere with the switching frequency tracking. It is used to guarantee the positive-definiteness of the Hessian, shown explicitly in Appendix~\ref{app:sphdec}. {After solving \eqref{eq:c2-opt} for its optimal solution $\b{U}_\augp^\opt(k) = [\b{u}_\augp(k)^\top\ \b{u}_\augp(k+1)^\top\ \ldots\ \b{u}_\augp(k+N_\p-1)^\top]^\top$, the \ac{c2} applies then the first input $\b{u}_\augp^\opt(k)$.} \looseness=-1

To solve the \ac{c2} in real-time, a sphere decoder similar to~\cite{geyer_2014, liegmann_2017} can be deployed.  The sphere decoding algorithm is a branch and bound technique applicable to \ac{ILS} problems\footnote{The \ac{ILS} problem is known to be NP-hard~\cite{van1981another}.} and is far more efficient than an exhaustive enumeration. The \ac{ILS} reformulation of the \ac{c2} and the algorithm are relegated to Appendix~\ref{app:sphdec}. Note that these were previously not derived or provided in~\cite{stellato_2017}. \looseness=-1

\section{FCS-MPC for $f_\sw$ Limiting via Slack Variables}
\label{sec:aug_fcs-mpc}
This section introduces \ac{c3}. We introduce a switching frequency constraint with the slack variable \looseness=-1
\begin{equation}\label{eq:slack_var}
s(\ell) =  \max\{f_\sw(\ell)- f_\sw^*,0\},
\end{equation} 
where $f_\sw^*$ is the upperbound on the switching frequency,
and penalize the slack variable in the cost function. We then develop a sphere decoder for solving the underlying problem and discuss a significant computational speed-up by utilizing lower bounds on this new cost. The \ac{c3} block diagram is identical to that of \ac{c2} depicted in Figure \ref{fig:MPC-tracking_block-diag}.\looseness=-1

\subsection{Problem formulation}
Similar to the $\b{p}(k)$ variable, the slack variable $s(k)$ also has a nonlinear output relation to the switch positions, which can be inferred from the definition in \eqref{eq:slack_var}. The nonlinearity originates from both the $\mathrm{max}$ operator and also the absolute value operator that is used to obtain the switching transition variable $\b{p}(k)$ appearing in \eqref{eq:5}. Thus, it would later not be possible to obtain an ILS reformulation unless the slack variable is included in the input vector. 

To this end, we augment the input vector as $\b{u}_\augs(k) = [\b{u}_\ph(k)^\top\ s(k+1)]^\top$. Given $\b{x}_\augs(k) = \b{x}_\ph(k)$ and $\b{y}_\augs(k) = \b{y}_\ph(k)$, the new dynamics are $\b{x}_\augs(k+1) = \b{A}_\augs\b{x}_\augs(k) + \b{B}_\augs\b{u}_\augs(k) + \b{D}_\augs\b{v}_\g(k)$ and $\b{y}_\augs(k) = \b{C}_\augs\b{x}_\augs(k)$.
Here, $\b{A}_\augs$, $\b{B}_\augs$, $\b{C}_\augs$, and $\b{D}_\augs$ can be defined simply by the physical system matrices in \eqref{eq:0}. \looseness=-1

We formulate the \ac{c3}'s optimization problem as \looseness=-1
\begin{equation}\label{eq:6}
\begin{split}
&\min \sum_{\ell=k}^{k+N_\p-1} || \b{y}_\augs^*(\ell+1) - \b{y}_\augs(\ell+1) ||_{\b{Q}_\augs}^2 \\
&\hspace{2.95cm} + \lambda_\sw ||s(\ell+1)||_2^2 + \lambda_\U || \Delta \b{u}_\ph(\ell) ||_2^2 \\
    &\ \mathrm {s.t.} \ \b{x}_\augs(\ell+1) = \b{A}_\augs\b{x}_\augs(\ell) + \b{B}_\augs\b{u}_\augs(\ell) + \b{D}_\augs\b{v}_g(\ell),\forall \ell, \\
    &\quad\quad\ \b{y}_\augs(\ell+1) = \b{C}_\augs \b{x}_\augs(\ell+1),\forall \ell, \\
    &\quad\quad\ \b{x}_\sw(\ell+1) = \b{A}_\sw\b{x}_\sw(\ell) + \b{B}_\sw \b{p}(\ell) ,\forall \ell, \\
    &\quad\quad\ s(\ell+1) = \max\left\{\b{C}_\sw\b{x}_\sw(\ell+1) - f_\sw^*,\ 0\right\},\forall \ell,\\
    &\quad\quad\ \b{u}_\ph(\ell) \in \{-1,\ 0,\ 1\}^3,\forall \ell, 
\end{split}
\end{equation} 
where $\b{Q}_\augs = \I_{2}$ and $\b{I}_n \in \mathbb{R}^{n\times n}$ denotes the identity matrix. Observe that the switching frequency penalty now relates to the penalization of the slack variable. 

We reformulate \eqref{eq:6} as an \ac{ILS} problem.  The full-horizon vector $\b{Y}_\augs(k+1)$ is a function of $\b{U}_\augs(k)$, $\b{x}_\augs(k)$, and $\b{v}_\g(k)$, i.e., $\b{Y}_\augs(k+1) = \b{\Gamma}_\augs \b{x}_\augs(k) + \b{\Upsilon}_\augs \b{U}_\augs(k) + \b{\Psi}_\augs \b{v}_\g(k).$ 

Other full-horizon variables $\b{S}(k+1)$ and $\Delta \b{U}_\augs(k)$ are defined as 
$\Delta \b{U}_\augs(k) = \b{\Pi}_\augs \b{U}_\augs(k) - \b{E}_\augs \b{u}_\augs(k-1)$, and $\b{S}(k+1) = \b{L}_\augs \b{U}_\augs(k)$, where $\b{\Gamma}_\augs$, $\b{\Upsilon}_\augs$, $\b{\Psi}_\augs$, $\b{\Pi}_\augs$, $\b{E}_\augs$, $\b{L}_\augs$, and $\bar{\b{Q}}_\augs$ are defined in Appendix~\ref{app:matrices} following the procedure originally described in \cite{geyer_2014}.

The objective in~\eqref{eq:6} is equivalent to
\begin{align*}
    J(\b{U}_\augs(k)) &= ||\b{Y}_\augs^\ast(k+1) - \b{Y}_\augs(k+1) ||_{\b{\bar{Q}}_\augs}^2 \\
    &\quad\quad + \lambda_\sw ||\b{S}(k+1)|| + \lambda_\U ||\Delta \b{U}_\augs(k)||,
\end{align*}
and can be regrouped as
\begin{align*}
    J(\b{U}_\augs(k)) = \b{U}_\augs(k)^\top \b{H}_\augs \b{U}_\augs(k) + 2\b{\Theta}_\augs(k)\b{U}_\augs(k) + \theta_\augs(k),
\end{align*}
where $\b{H}_\augs = \b{\Upsilon}_\augs^\top \bar{\b{Q}}_\augs \b{\Upsilon}_\augs + \lambda_\sw \b{L}_\augs^\top \b{L}_\augs + \lambda_\U \b{\Pi}_\augs^\top\b{\Pi}_\augs$ is called the Hessian and $\b{\Theta}_\augs(k) = ((\b{\Gamma}_\augs \b{x}_\augs(k)\allowdisplaybreaks - \b{Y}_\augs^*(k+1)\allowdisplaybreaks + \b{\Psi}_\augs \b{v}_\g(k))^\top \bar{\b{Q}}_\augs\b{\Upsilon}_\augs\allowdisplaybreaks- \lambda_\U(\b{E}_\augs\b{u}_\augs(k-1))^\top \b{\Pi}_\augs )^\top)^\top$ is the linear part of the cost. The cost term $\theta_\augs(k)$ is independent of the decision variable $\b{U}_\augs(k)$ and, therefore, gets discarded, and instead gets replaced by another constant term to complete the squares. \looseness=-1

The matrices $\b{\Upsilon}_\augs$, $\b{\Pi}_\augs$, and $\b{L}_\augs$ are constructed so that $\lambda_\sw \b{L}_\augs^\top \b{L}_\augs + \lambda_\U \b{\Pi}_\augs^\top\b{\Pi}_\augs \succ 0$ and $\b{\Upsilon}_\augs^\top \bar{\b{Q}}_\augs \b{\Upsilon}_\augs \succeq 0$. Consequently, it holds for the Hessian that $\b{H}_\augs \succ 0$, and there exists a generator matrix $\b{V}_\augs$ such that $\b{V}_\augs^\top\b{V}_\augs = \b{H}_\augs$. The matrix $\b{V}_\augs$ is lower triangular, and every fourth row has only the slack variable {penalty} as a diagonal entry. 

With the definitions above, \eqref{eq:6} is equivalent to the following \ac{ILS} problem: \looseness=-1
\begin{equation}\label{eq:9}
\begin{split}
    &\min_{\b{U}_\augs(k)} ||\hat{\b{U}}_\augs(k) - \b{V}_\augs\b{U}_\augs(k)||_2^2 \\
    &\quad\mathrm{s.t.}\ \b{u}_\ph(\ell) \in \{-1,\ 0,\ 1\}^{3}, \forall \ell,
    \\ 
    &\quad\quad\quad s(\ell) = \tilde{h}_{\ell}(\b{U}_{\past}(\ell)), \forall \ell>k.
\end{split}
\end{equation}
where $\hat{\b{U}}_\augs(k) = -\b{V}_\augs \b{H}_{\augs}^{-1} \b{\Theta}_{\augs}(k)$ is the unconstrained solution to~\eqref{eq:6}, i.e., the solution of~\eqref{eq:6} without integer and slack constraints. 
To simplify the presentation, the slack variables' functional dependence on the inputs is captured by the functions $\tilde{h}_\ell:\mathbb{R}^{3(\ell-k)} \mapsto \mathbb{R}_{\geq 0}$, resulting in the constraint $s(\ell) = \tilde{h}_\ell(\b{U}_{\past}(\ell)), \forall \ell>k$, where $\b{U}_{\past}(\ell) = [\b{u}_\ph(k)^\top \ldots \b{u}_\ph(\ell-1)^\top]^\top, \forall \ell>k$ denotes the past inputs.\footnote{This change in notation highlights that the sphere decoder in the next subsection will in fact be applicable to slack variables of constraints with any nonlinear relationship to the past inputs $\b{U}_{\past}(\ell)$.} These functions depend on the initial states, e.g., $\b{x}_\sw(k)$, and can be evaluated iteratively via the constraints in~\eqref{eq:6} for different~$\ell$.
For the slack variables $s(\ell)$, corresponding entries of the unconstrained solution $\hat{{\b{U}}}_\augs$ are 0, since the constraint $s(\ell) = \tilde{h}_\ell(\b{U}_{\past}(\ell))$ is not imposed on the unconstrained solution.
\looseness=-1

\subsection{Sphere decoder}
The objective in \eqref{eq:9} relates to the squared distance between the unconstrained solution $\hat{\b{U}}_\augs(k)$ and the optimization variable $\b{U}_\augs(k)$ in the space of the lattice defined by $\b{V}_\augs$. Given a feasible initial solution $\b{U}_{\augs,\ini}(k)$ with the objective value $d^2_{\ini}$, the sphere decoder systematically searches for a better feasible solution within a sphere with the initial squared radius of $\rho^2=d^2_{\ini}$. The sphere decoder shrinks the squared radius $\rho^2$ whenever it finds a better solution $\b{U}_\augs(k)$. Through this process, suboptimal solutions are discarded until only the optimal solution $\b{U}_\augs^\opt(k)$ is left inside the sphere. \looseness=-1

% The sphere decoder exploits the lower triangular property of the generator matrix $\b{V}_\augs$ to compute the squared distance $d^2$ \textcolor{blue}{of a solution candidate $\b{U}_\augs(k)$} in the form % 
The exploration procedure of the sphere decoder for a better solution relies on the structure of the squared distance~$d^2$. Let $\b{\Xi}_{i}$ be the $i^{\th}$ row vector of a matrix $\b{\Xi}$, and $\varphi_{i}$ be the $i^{\th}$ element of a vector $\b{\varphi}$, then \looseness=-1
\begin{equation*}\label{eq:d^2}
    d^2 =  \sum\limits_{i=1}^{N_\p}\left|\left| \hat{U}_{\augs, i} - \b{V}_{\augs, i}\b{U}_\augs(k) \right|\right|_2^2.
\end{equation*}
Instead of directly calculating $d^2$ for the full candidate solution vector, we can iteratively compute an incremental squared distance term $d_l^2$, where $l\in\{1,\ldots,4 N_\p\}$ and $d_{4N_\p}^2 = d^2$. For this, we exploit the lower triangular structure of $\b{V}_\augs$, i.e., \looseness=-1
\begin{align}
    d_l^2 = \sum\limits_{i=1}^{l}
    \left(\hat{U}_{\augs, i} - 
    \sum\limits_{j=1}^{i} V_{\augs, i,j}U_{\augs, j}(k)\right)^2, \label{eq:4}
\end{align}{where $\Xi_{i,j}$ denotes the element in the $i^{\th}$ row and $j^{\th}$ column of a matrix $\b{\Xi}$. The incremental squared distance $d_l^2$ allows the sphere decoder to avoid evaluating the objective function at once for a fully-determined candidate solution. Instead, we can determine the entries of the candidate vector one by one, starting at time step $k$. We can then increment the squared distance accordingly. This way, a large set of candidates can be discarded whenever the incremental distance already exceeds the sphere's current radius, $d_l^2 > \rho^2$, before all the entries are determined.

% \textcolor{blue}{The index $l$ can also be interpreted as a node level in a search tree of all possible solution candidates~\cite{dorfling_2020}. In the tree, every first three levels correspond to the switch positions $\b{u}_\ph(\ell) \in \{-1,0,1\}^{3}$ and every fourth level to the slack variable $s(\ell+1) \in \mathbb{R}_{\geq 0}$.} \looseness=-1

With this in mind, the sphere decoder can be formulated as described in Algorithm \ref{alg:1}. When starting the recursive algorithm, we let $\b{U}_\augs^\opt = \b{U}_{\augs}^\ini(k)$, $\hat{\b{U}}_\augs = \hat{\b{U}}_\augs(k)$, $d_l^2 = 0$, $l=1$, $\ell=k$, $\b{u}_{\augs}^{\prev} = \b{u}_\augs(k-1)$, and $\rho^2 = \rho_\ini^2(k)$. 
The initial solution $\b{U}_{\augs}^\ini(k)$ can be obtained with an educated guess (i.e., a shifted version of the previous optimal solution $\b{U}_\augs^\opt(k-1)$), the Babai estimate (e.g., see~\cite{dorfling_2020}), or the best of both. Additionally, we need to pre-compute the corresponding slack variables in $\b{U}_\augs^\ini(k)$. \looseness=-1

\begin{algorithm}[t!]
\caption{Modified sphere decoder for \ac{c3}}
\label{alg:1}
{\setstretch{.99}\begin{algorithmic}[1]
    \State Function: FLSphDec
    \State Input: $\b{U}_\augs$, $\b{U}_\augs^\opt$, $\hat{\b{U}}_\augs$, $\rho^2$, $d_{l-1}^2$, $\underbar{d}_{\ell}^2$, $l$, $\ell$, $\b{u}_{\ph}^{\prev}$, $\b{x}_\sw(\ell)$
    \State Output: $\b{U}_\augs^\opt$, $\rho^2$
    \If{$\mod(l,4)\neq0$}
    \State $\mathcal{U} \gets \{-1,0,1\}$
    \Else
    \State $\b{u}_{\ph}^{\prev} \gets \b{u}_\ph(\ell)$
    \State $\ell \gets \ell + 1$
    \State Compute $\b{x}_\sw(\ell)$, $s(\ell)$ and $\underbar{d}_\ell^2$, and let $\mathcal{U} \gets s(\ell)$ \label{line:slack}
    \EndIf
    \ForAll{$u\in\mathcal{U}$}
    \State $U_{\augs,l} \gets u$
    % \State $d_{l}^2 \gets || \hat{U}_{\augs,l} - \b{V}_{\augs,l}\b{U}_\augs ||_2^2 + d_{l-1}^2$
    \State $d_{l}^2 \gets d_{l-1}^2 + \left(\hat{U}_{\augs,l} - \sum\nolimits_{j=1}^{l} V_{\augs, l,j}U_{\augs, j}(k)\right)^2$
    \If{$d_{l}^2 + \underbar{d}_{\ell}^2 < \rho^2$}
    \If{$l < 4N_\p$}
    % \State Update $\b{u}_{\ph}^{\prev}$ if mod$(l,4) = 0$
    \State New recursion: $[\b{U}_\augs^\opt$, $\rho^2]$ $\gets$ FLSphDec( $\b{U}_\augs$, $\b{U}_\augs^\opt$, $\hat{\b{U}}_\augs$, $\rho^2$, $d_l^2$, $\underbar{d}_{\ell}^2$, $l+1$, $\ell$, $\b{u}_{\ph}^{\prev}$, $\b{x}_\sw(\ell)$)
    \Else
    \State $\b{U}_\augs^\opt \gets \b{U}_\augs$ and $\rho^2 \gets d_l^2$
    \EndIf
    \EndIf
    \EndFor
\end{algorithmic}}
\end{algorithm}

A tree of depth $4N_\p$, and levels denoted by $l$, representing all possible solutions of $\b{U}_\augs(k)$, is traversed. Every first three nodes, i.e., $\mod(l,4)\neq0$, decide upon the switch position $\b{u}_\ph(\ell)$, whereas every fourth node computes the corresponding slack variable $s(\ell+1)$. This can be done by evaluating functions $\tilde{h}_\ell$. For~\eqref{eq:9}, this is done with the definition of the slack variable in \eqref{eq:slack_var} after iteratively computing $\b{x}_\sw(\ell+1)$ via~\eqref{eq:5}. 
Note that the slack variable prediction models are not explicitly included when deriving~$\b{V}_\augs$, but remain hidden in the computation in Line~\ref{line:slack}. A simple lower bound for the slack variable cost to be incurred in the future steps within the prediction horizon is given by $\underbar{d}^2_\ell=0$. In the next subsection, to cut branches as early as possible, we will provide a tighter lower bound $\underbar{d}_{\text{bound},\ell}^2$. This bound will be derived based on the prediction model of the slack variable for the switching frequency constraint. The algorithm can be used with either of $\underbar{d}^2_\ell=0$ or $\underbar{d}^2_\ell=\underbar{d}_{\text{bound},\ell}^2$.
\looseness=-1

After picking the value for the input, the incremental cost~\eqref{eq:4} at level $l$ is compared with the sphere's squared radius $\rho^2$, i.e., the current best solution. The tree is further traversed if $d_l^2+\underbar{d}_l^2 < \rho^2$ still holds. Otherwise, the node with its sub-tree is pruned. \looseness=-1

This procedure is repeated until either all branch values of a node are cut or a solution with $d_{4N_\p}^2 < \rho^2$ is found at a leaf node. If $d_{4N_\p}^2 < \rho^2$, then the sphere radius is reduced to $\rho^2 = d_{4N_\p}^2$ and $\b{U}_\augs^\opt$ is updated. The algorithm starts back-tracking. If a new solution with cost $(d_i^\prime)^2 < \rho^2$ is found, the algorithm traverses the new branch again. The optimal solution $\b{U}_\augs^\opt(k) = \b{U}_\augs^\opt$ is found once the last branch has been cut and the node at $l=1$ is reached.

%Next, we will provide a tighter lower bound $\underbar{d}_{\text{bound},\ell}^2$ to cut branches as early as possible. Algorithm~\ref{alg:1} can be initialized with $\underbar{d}_\ell^2 = 0$ or with the tighter bound $\underbar{d}_\ell^2 = \underbar{d}_{\text{bound},k}^2$.}

\subsection{Lower bound for computational speed-up}
Consider $\ell<k+N_\p$, since otherwise there is no future cost to be incurred by the slack variable. Let \begin{equation}\label{eq:slack_cost}
    J_\augs(\ell) = \lambda_\sw\allowdisplaybreaks \sum_{n=1}^{k+N_\p-\ell} \allowdisplaybreaks ||\max(\b{C}_\sw \b{x}_\sw(\ell + n) - f_\sw^\ast, 0) ||_2^2,
\end{equation} denote the cost incurred from the slack variable (strictly) after time step~$\ell$ if generic inputs $\b{u}_\ph(\ell-1+ n)$, $n\in\{1,\ldots, k+N_\p-\ell\}$, are used for the rest of the horizon. We show that there exists a lower bound $\underbar{d}_{\text{bound},\ell}^2 \leq J_\augs(\ell)$, which is often nontrivial, i.e., $\underbar{d}_{\text{bound},\ell}^2>0$, and can simply be computed from the already determined state variable $\b{x}_\sw(\ell)$ in step $\ell$.
Intuitively, this will originate from the case where no more switching occurs in the remainder of the prediction horizon of the algorithm. 

\begin{proposition}\label{prop:lwb}
Given $\ell<k+N_\p$, define $$\underbar{d}_{\text{bound},\ell}^2= \lambda_\sw \sum\limits_{n=1}^{k+N_\p-\ell} ||\max(\b{C}_\sw \b{A}_\sw^{n} \b{x}_\sw(\ell)-f_\sw^*, 0)||_2^2.$$ We have that $\underbar{d}_{\text{bound},\ell}^2\leq J_\augs(\ell).$
\end{proposition}
\begin{IEEEproof}
Observe that we obtain the switching state $\b{x}_\sw(\ell+n)$ in a sequential manner, i.e., $$\b{x}_\sw(\ell+n) = \b{A}_\sw^{n} \b{x}_\sw(\ell) + \sum_{j=0}^{n-1} \b{A}_\sw^{j}\b{B}_\sw |\Delta \b{u}_\ph(\ell+j)|.$$ We have that $$\b{A}_\sw^n \b{x}_\sw(\ell) \leq \b{A}_\sw^n \b{x}_\sw(\ell) + \sum_{j=0}^{n-1} \b{A}_\sw^j \b{B}_\sw |\Delta \b{u}_\ph(\ell+j)|,$$ due to non-negativity of the term on the right.
Moreover, $\max((\cdot) - f_\sw^\ast, 0)$ is a non-negative and a non-decreasing operator. 

Invoking these two observations, the bound $\underbar{d}_{\text{bound},\ell}^2$ below obtained by assuming no future switchings can be shown to provide a lower bound to \eqref{eq:slack_cost}, that is, the slack variable cost to be incurred in the future steps for any generic future inputs: \looseness=-1
\begin{equation*}\label{eq:13}
\begin{split}
    \underbar{d}_{\text{bound},\ell}^2 &= \lambda_\sw \sum\limits_{n=1}^{k+N_\p-\ell} ||\max(\b{C}_\sw \b{A}_\sw^{n} \b{x}_\sw(\ell)-f_\sw^*, 0)||_2^2\\
    &\leq \lambda_\sw\allowdisplaybreaks \sum_{n=1}^{k+N_\p-\ell} \allowdisplaybreaks ||\max(\b{C}_\sw \b{x}_\sw(\ell + n) - f_\sw^\ast, 0) ||_2^2\\
    &=J_\augs(\ell).
    \end{split}
\end{equation*}
This concludes the proof.
\end{IEEEproof}

{\setstretch{.99}\begin{table}[t!]
     \centering
     \caption{Simulation Parameters}
     \label{tab:0}
     {\setstretch{.99}{\footnotesize\begin{tabular}{cccccc}
        \hline
        \multicolumn{2}{c}{Physical System} & \multicolumn{2}{c}{Controllers} & \multicolumn{2}{c}{Simulation}\\
        \hline
        $L$ & \SI{.266}{\pu} & $N_\p$ & 5 & $T_\Sim$ & \SI{.5}{\micro\second} \\
        $R$ & \SI{.015}{\pu} & $a_1$ & $0.99$ & $a_1^\text{visual}$ & $0.995$ \\
        $V_\g$ & \SI{1}{\pu} & $a_2$ & $0.99$ & $a_2^\text{visual}$ & $0.995$ \\
        $V_\dc$ & \SI{1.9}{\pu} & $\lambda_\U$ & $13e^{-3}$ &  &  \\
        $f_1$ & \SI{50}{\hertz} & $\lambda_\sw$ & $60$ &  &  \\
         &  & $T_\C$ & $\SI{100}{\micro\second}$ &  & \\
        \hline
     \end{tabular}}} \vspace{.1cm}
 \end{table}}
	
\section{Case Studies}\label{sec:numerics}
We benchmark the \ac{c3} against \ac{c2}. We showcase the benefits of \ac{c3} in common test cases such as steady-state current tracking and power ramp-ups. We also demonstrate the controllers during changes in the switching frequency reference to showcase the impact of the different objective functions. \looseness=-1

For the system presented in Section~\ref{sec:system_modeling}, the parameters are shown in Table~\ref{tab:0}. Here, $T_\Sim$ refers to the simulation step size, $a_1^\text{visual}$ and $a_2^\text{visual}$ are used for plotting purposes. The controller's current tracking performance is evaluated with the current \ac{TDD}. The current \ac{TDD} is well-suited for steady-state measurements, however, it is less applicable for evaluating the current tracking performance over short intervals during transients. Hence, we define the $2$-norm of the current tracking error: 
$$e_I^2 = \frac{1}{\sqrt{2} T} \sum\limits_{t=t_0}^{t_0+T}||\b{x}_\ph^*(t) - \b{x}_\ph(t)||_2^2.$$ The average switching frequency $\bar{f}_\sw$ is defined as $$\bar{f}_\sw = \frac{1}{T}\sum\limits_{t=t_0}^{t_0 + T} f_\sw(t).$$ The two controllers are tuned with exactly the same parameters for a fair comparison. \looseness=-1

\begin{table}[t!]
	\centering
	\caption{Simulation results; Current TDD and average switching frequency $\bar{f}_\sw$ for the steady-state simulation.}
	\label{tab:simulation}
	{\setstretch{.99}{\footnotesize\begin{tabular}{cccc}
        \hline
        & \ac{c2} & \ac{c3} & Improvement [\%]\\
        \hline
        Current TDD [\%] & 4.95 & 4.70 & 5.1 \\
        Avg. sw. frequency [Hz] & 253 & 248 & - \\
		\hline
	\end{tabular}}} \vspace{.1cm}
\end{table}

\begin{table}[t!]
	\centering
	\caption{Simulation results; Errors $e_{I,\ct}$ and $e_{I,\cs}$ for the power ramp and the switching frequency step-up.}
	\label{tab:simulation}
	{\setstretch{.99}{\footnotesize\begin{tabular}{cccc}
        \hline
		& $e_{I,\ct}\ [\%]$ & $e_{I,\cs}\ [\%]$ & Improvement [\%] \\
		\hline
		Power ramp-up & $6.8$ & $6.1$ & $9.8$ \\
        $f_\sw$ step & $8.8$ & $3.9$ & $56$ \\
		\hline
	\end{tabular}}}
 \vspace{.1cm}
\end{table}

\begin{figure}[t!]
	\centering
	\begin{subfigure}{.49\columnwidth}
		\centering
		\setlength{\fwidth}{1.1\linewidth}
		\setlength{\fheight}{0.9\fwidth}
		\input{Figures/f_sw-ss.tex}
  \vspace{-.6cm}
		\caption{Switching frequencies}
		\label{fig:f_sw-ss}
	\end{subfigure}
	% \hspace{.1cm}
	\begin{subfigure}{.49\columnwidth}
        \centering
    	\setlength{\fwidth}{1.1\linewidth}
    	\setlength{\fheight}{.9\fwidth}
        \include{Figures/i-ss.tex}
        \vspace{-1cm}
        \caption{Currents over a fundamental}
        \label{fig:i-ss}
    \end{subfigure}
    \caption{Steady-state operation.}
\end{figure}
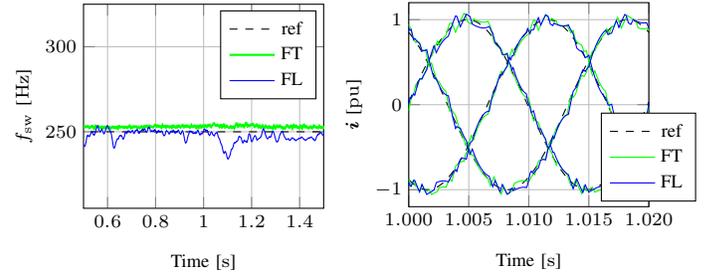

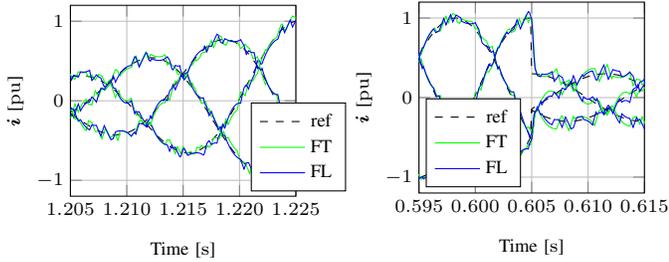
\begin{figure}[t!]
	\centering
	\begin{subfigure}{.47\columnwidth}
		\centering
		\setlength{\fwidth}{1.1\linewidth}
		\setlength{\fheight}{0.9\fwidth}
		\input{Figures/i-iramp.tex}
  \vspace{-.2cm}
		\caption{Currents during power ramp.}
		\label{fig:i-ramp}
	\end{subfigure}
	\quad
	\begin{subfigure}{.47\columnwidth}
        \centering
    	\setlength{\fwidth}{1.1\linewidth}
    	\setlength{\fheight}{.9\fwidth}
        \include{Figures/i-istep}
\vspace{-1cm}
        \caption{Currents during power step.}
        \label{fig:i-step}
    \end{subfigure}
    \caption{Current transients.}
\end{figure}

We first perform a steady-state simulation for \SI{1.5}{\second} (i.e., 75 fundamental periods). An outer controller generates a current reference, injecting to the grid $P=\SI{1}{\pu}$, $Q=\SI{0}{\pu}$. We initialize the measurement at $t_0 = \SI{.5}{\second}$ and consider $T = \SI{1}{\second}$. In the steady-state, we measure an average switching frequency of $\bar{f}_{\sw,\ct} = \SI{253}{\hertz}$ for \ac{c2} and $\bar{f}_{\sw,\cs} = \SI{248}{\hertz}$ for \ac{c3}, with a reference of $f^*_\sw = \SI{250}{\hertz}$. The difference in average switching frequencies is due to the different objectives and it is visualized in Figure~\ref{fig:f_sw-ss}. The current \ac{TDD} measurement is provided in Table~\ref{tab:simulation} and shows a $5.1\%$ improvement for the \ac{c3} over the \ac{c2} at a lower average switching frequency. Figure~\ref{fig:i-ss} shows the currents over one fundamental period during steady-state. \looseness=-1

Similar improvements can also be observed when applying an active power ramp. Let $t_0 = \SI{1.205}{\second}$, and consider one fundamental period $T_1 = \frac{1}{f_1}$ length. The ramp is initialized at exactly $t_0 = \SI{1.205}{\second}$ with $P = \SI{.3}{\pu}$ and ends after one fundamental period at $P = \SI{1}{\pu}$, while $Q = \SI{0}{\pu}$. The simulation results with the improved current tracking performance of the \ac{c3} over \ac{c2} are shown in Figure~\ref{fig:i-ramp} and Table~\ref{tab:simulation}. Additionally, Figure~\ref{fig:i-step} shows a power step, in which both controllers perform equally well. These different scenarios show that \ac{c3} is at least as good, and if not better, in tracking a current reference when compared to \ac{c2}. This is not surprising, considering that \ac{c2} might unnecessarily prioritize switching frequency tracking. \looseness=-1

Next, though it is not realistic, to showcase an obvious case of current tracking improvement, we apply a step in the switching frequency reference. The simulation results are depicted in Figure~\ref{fig:f_sw-step}. Here, the \ac{c2} realizes a faster step than \ac{c3} in terms of switching frequency. We initialize the measurement at $t_0 = \SI{.4}{\second}$ and measure one fundamental period. Table~\ref{tab:simulation} shows a significantly improved performance for \ac{c3} because \ac{c2} jeopardizes current tracking to increase the switching frequency as fast as possible.\looseness=-1 

To evaluate the benefits of the computational speed-up, we simulated the system while measuring the total time for sphere decoding.\footnote{All problems are solved via MATLAB on a computer equipped with \SI{16}{\giga\byte} RAM and a \SI{1.8}{\giga\hertz} quad-core Intel i7 processor.} The results are shown in Table~\ref{tab:speed-up}. Here, we can see a clear improvement in the \ac{c3} computation times; long computation instances are eliminated thanks to the lower bound in Proposition~\ref{prop:lwb}. These reduced computation times are comparable to those of \ac{c2}. A final comparison of the characteristics of the two controllers is provided in Table~\ref{tab:final}. \looseness=-1

{\setstretch{.99}\begin{table}[t!]
	\centering
	\caption{Total simulation time (Sim time), maximum solving time (Max time), $70^{\text{th}}$ percentile, and $95^{\text{th}}$ percentile of \ac{c2}, and \ac{c3} with and without speed-up.}
	\label{tab:speed-up}
	{\setstretch{.97}{\footnotesize\begin{tabular}{cccc}
        \hline
		& \ac{c2} & \ac{c3} no spd.-up & \ac{c3} w/ spd.-up\\
		\hline
		Sim time [s] & $30$ & $78$ & $22$ \\
		Max time [ms] & $9.6$ & $232$ & $7.7$ \\
		$70^{\text{th}}$ per. [ms] & $0.60$ & $0.87$ & $0.55$ \\
		$95^{\text{th}}$ per. [ms] & $1.40$ & $13.50$ & $1.40$ \\
		\hline
	\end{tabular}}} \vspace{.1cm}
\end{table}}

\begin{figure}[t!]
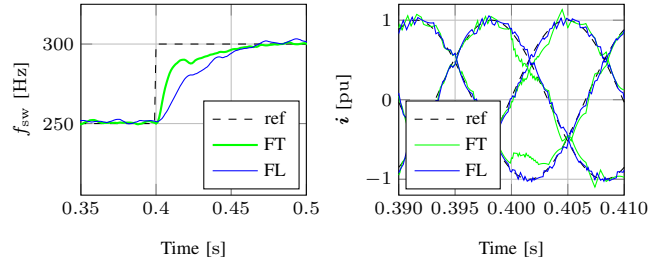

	\centering
	\begin{subfigure}{.47\columnwidth}
        \centering
    	\setlength{\fwidth}{1.1\linewidth}
    	\setlength{\fheight}{.9\fwidth}
        \include{Figures/fsw-fswstep}
        \vspace{-1cm}
        \caption{Switching frequencies}
        \label{fig:f_sw-fswstep}
	\end{subfigure}
	% \quad
	\begin{subfigure}{.47\columnwidth}
        \centering
    	\setlength{\fwidth}{1.1\linewidth}
    	\setlength{\fheight}{.9\fwidth}
        \include{Figures/i-fswstep}
        \vspace{-1cm}
        \caption{Currents}
        \label{fig:i-fswstep}
    \end{subfigure}
    \caption{Switching frequency step-up.}
    \label{fig:f_sw-step}
\end{figure}

{\begin{table}[t!]
	\centering
	\caption{Comparison of the characteristics of controllers.}
	\label{tab:final}
	{{\footnotesize\begin{tabular}{ccc}
		\hline  \vspace{.1cm}
		& \ac{c2} & \ac{c3}  \\ 
		\hline
  \vspace{.1cm}
		Cost function & \begin{tabular}{@{}c@{}} Current tracking \\ $f_\sw$ tracking \end{tabular}& \begin{tabular}{@{}c@{}} Current tracking \\ $f_\sw$ limiting \end{tabular}  \\
   \vspace{.1cm}
        Weighting factors & \begin{tabular}{@{}c@{}} $\lambda_\sw$\text{: $f_\sw$ tracking} \\  \quad weight \\ $\lambda_\U$\text{: switching penalty}\end{tabular} & \begin{tabular}{@{}c@{}} $\lambda_\sw$\text{: $f_\sw$ constraint} \\ \quad\ slack weight\\$\lambda_\U$\text{: switching penalty}\end{tabular}  \\  \vspace{.1cm}
		\begin{tabular}{@{}c@{}} ILS optimizer \\ dimension\end{tabular} & $\mathrm{dim}(\b{U}_\augp(k))=6N_\p$ & $\mathrm{dim}(\b{U}_\augs(k))=4N_\p$  \\  \vspace{.1cm}
		Computational time & Low & \begin{tabular}{@{}c@{}} Similar to \ac{c2} \\ via speed-up (Prop.~\ref{prop:lwb})\end{tabular}\\  \vspace{.1cm}
		\begin{tabular}{@{}c@{}} Transient current \\ tracking performance\end{tabular} & Good & Improved  \\
		\hline
	\end{tabular}}}
\end{table}}

\vspace{-0.1cm}
\section{Conclusion}
We introduced \ac{c3}, replacing the switching frequency tracking objective of \ac{c2} with switching frequency limiting. \ac{c3} has been shown to perform better than \ac{c2} in multiple scenarios. Additionally, a computational speed-up was proposed for the sphere decoder of \ac{c3}. Future work could implement \ac{c3} on an embedded FPGA.\looseness=-1 

\vspace{-0.1cm}	
\bibliographystyle{IEEEtran}
\bibliography{brief}

\appendices

\section{Sphere decoder for \ac{c2}}\label{app:sphdec}
Besides the conventional FCS-MPC of \cite{geyer_2014}, the \ac{c2} optimization problem is also amenable to an \ac{ILS} reformulation. We write the full-horizon vector $\b{Y}_\augp(k+1)$ as a function of $\b{U}_\augp(k)$, $\b{x}_\augp(k)$, $\b{v}_\g(k)$, that is, 
$$\b{Y}_\augp(k+1) = \b{\Gamma}_\augp \b{x}_\augp(k) + \b{\Upsilon}_\augp \b{U}_\augp(k) + \b{\Phi}_\augp \b{v}_\g(k).$$ 
We also define the vector $\b{P}(k) = [\Delta \b{u}_\ph^\top(k)\ \b{p}^\top(k) \ldots \Delta \b{u}_\ph^\top(k+N_\p-1)\ \b{p}^\top(k+N_\p-1)]^\top$, which can be computed as $\b{P}(k) = \b{\Pi}_\augp\b{U}_\augp(k) - \b{E}_\augp\b{u}_\augp(k-1).$
We remark that this definition slightly deviates from the original notation we have for the full-horizon vectors.
The matrices $\b{E}_\augp$, $\b{\Pi}_\augp$, $\b{\Gamma}_\augp$, $\b{\Upsilon}_\augp$, $\b{\Phi}_\augp$, $\bar{\b{Q}}_\augp$ are defined in Appendix \ref{app:matrices}. \looseness=-1

{We use the full-horizon vectors to write the objective~\eqref{eq:c2-opt} as \looseness=-1
\begin{align*}
\begin{split}
    J(\b{U}_\augp(k)) &= ||\b{Y}_\augp^\ast(k+1) - \b{Y}_\augp(k+1) ||_{\b{\bar{Q}}_\augp}^2 \\
    &\quad\quad\quad\quad\quad\quad + \frac{\lambda_\U}{2} ||\b{P}(k+1)||_2^2.
\end{split}
\end{align*}
As a remark, $||\b{P}(k)||_2^2  = 2||\Delta \b{U}(k)||_2^2$ because $||\Delta \b{u}_\ph(\ell)||_2^2 = ||\b{p}(\ell)||_2^2$. A similar trick was utilized also  in~\cite{liegmann_2017}. By using the definitions of $\b{Y}_\augp(k+1)$ and $\b{P}(k)$, we obtain the quadratic form} \looseness=-1
\begin{equation*}
    J(\b{U}_\augp(k)) = \b{U}_\augp(k)^\top \b{H}_\augp \b{U}_\augp(k) + 2\b{\Theta}_\augp(k)\b{U}_\augp(k) + \theta_\augp(k),
\end{equation*}
where the Hessian is defined as $\b{H}_\augp = \b{\Upsilon}_\augp^\top\bar{\b{Q}}_\augp \b{\Upsilon}_\augp + \frac{\lambda_\U}{2} \b{\Pi}_\augp^\top\b{\Pi}_\augp$ and the linear cost term is $\b{\Theta}_\augp(k) = (\b{\Gamma}_\augp \b{x}_\augp(k) - \b{Y}_\augp^*(k+1) + \b{\Psi}_\augp \b{v}_\g(k))^\top \bar{\b{Q}}_\augp\b{\Upsilon}_\augp - \frac{\lambda_\U}{2} (\b{E}_\augp\b{u}_\augp(k-1))^\top \b{\Pi}_\augp )^\top$. The cost term $\theta_\augp(k)$ is independent of the decision variable $\b{U}_\augp(k)$ and, therefore, gets discarded. Instead, we replace it with another constant term to complete the squares. \looseness=-1

The \ac{ILS} reformulation of the \ac{c2} problem is\looseness=-1
\begin{equation}\label{eq:12}
\begin{split}
    &\min_{\b{U}_\augp(k)} ||\hat{\b{U}}_{\augp}(k) - \b{V}_\augp\b{U}_\augp(k)||_2^2 \\
    &\quad \mathrm{ s.t.}\ \b{u}_\ph(\ell) \in \{-1,\ 0,\ 1\}^{3},\forall \ell,\\
    &\quad\quad\quad \b{p}(\ell) = |\Delta \b{u}_\ph(\ell)|,\forall \ell,
\end{split}
\end{equation}
where the unconstrained solution $\hat{\b{U}}_\augp(k) = -\b{V}_\augp \b{H}^{-1} \b{\Theta}(k)$ is defined by the generator matrix $\b{V}_\augp$ from the Cholesky decomposition of $\b{H}_\augp = \b{V}_\augp^\top\b{V}_\augp$. For the absolute switching transitions $\b{p}(\ell)$, corresponding entries of the unconstrained solution $\hat{{\b{U}}}_\augp(k)$ are uncorrelated to the $\b{u}_\ph(\ell)$ entries, since the constraint $\b{p}(\ell) = |\Delta \b{u}_\ph(\ell)|$ is not imposed. \looseness=-1

\begin{algorithm}[t!]
    \caption{Modified sphere decoder for \ac{c2}}
    \label{alg:2}
    \begin{algorithmic}[1]
        \State Function: FTSphDec
        \State Input: $\b{U}_\augp$, $\b{U}_\augp^\opt$, $\hat{\b{U}}_\augp$, $\rho^2$, $d_{l-1}^2$, $l$, $\ell$, $k$, $\b{u}_{\ph}^{\prev}$
        \State Output: $\b{U}_\augp^\opt$, $\rho^2$
        \State $j=l-6(\ell-k)$
        \State $\mathcal{U} \gets \left\{ \begin{array}{l}
            \{-1,0,1\}, \text{ if } j \in \{1,2,3\}\\
            \{|u_{ j-3}(\ell) - u_{ j-3}^\prev|\}, \text{ otherwise}
        \end{array} \right.$
        \If{$j=6$} 
        \State $\b{u}_\ph^\prev \gets \b{u}_\ph(\ell)$, $\ell \gets \ell+1$
        \EndIf
        \ForAll{$u\in\mathcal{U}$}
        \State $U_{\augp, l} \gets u$
        % \State $d_{l}^2 \gets || \hat{U}_{\augp, l} - \b{V}_{\augp, l}\b{U}_\augp ||_2^2 + d_{l-1}^2$
        \State $d_{l}^2 \gets d_{l-1}^2 + \left(\hat{U}_{\augp,l} - \sum\nolimits_{j=1}^{l} V_{\augp, l,j}U_{\augp, j}(k)\right)^2$
        \If{$d_{l}^2 < \rho^2$}
        \If{$l < 6N_\p$}
        \State New recursion: $[\b{U}_\augp^\opt$, $\rho^2]$ $\gets$ FTSphDec( $\b{U}_\augp$, $\b{U}_\augp^\opt$, $\hat{\b{U}}_\augp$, $\rho^2$, $d_{l}^2$, $l+1$, $\ell$, $k$, $\b{u}_\augp^\prev$)
        \Else
        \State $\b{U}_\augp^\opt \gets \b{U}_\augp$ and $\rho^2 \gets d_l^2$
        \EndIf
        \EndIf
        \EndFor
    \end{algorithmic}
\end{algorithm}
To solve the \ac{ILS} problem presented in \eqref{eq:12}, a modified version of the sphere decoder similar to that from~\cite{geyer_2016,liegmann_2017} can be used. It is provided in Algorithm \ref{alg:2} for the sake of completeness. The algorithm iterates over inputs, starting from $\ell=k$, exploiting the lower triangular structure of the generator matrix  $\b{V}_\augp$. Here, $\b{V}_{\augp, l}$ refers to the $l^\mathrm{th}$ row of the matrix $\b{V}_{\augp}$, and similarly, $\hat{U}_{\augp, l}$ is the $l^\mathrm{th}$ entry of the vector $\hat{\b{U}}_{\augp}$. While maintaining the conventional algorithm structure, the input set $\mathcal{U}$ can simply be modified for the iterations involving the variable $\b{p}(k)$. The sphere decoder would first iterate over the input $u_{j}(\ell)$ for some $j\in\{a,b,c\}$ at a certain time step $\ell$. The absolute switching transition variables $p_j(\ell)$ would then simply be computed based on the switch positions $u_{j}(\ell)$ and $u_{j}(\ell-1)$. \looseness=-1

The algorithm is initialized with $\b{U}_\augp = \b{U}_\augp^\ini(k)$, $\hat{\b{U}}_\augp = \hat{\b{U}}_\augp(k)$, $\rho^2 = \rho_\ini^2(k)$, $d_{l-1}^2 = 0$, $l=1$, $\ell = k$, $\b{u}_\ph^\prev = \b{u}_\ph(k-1)$. The initial solution $\b{U}_\augp^\ini(k)$ can be computed with the educated guess (i.e., a shifted version of the previous optimal solution $\b{U}_\augp^\opt(k-1)$), Babai estimate (a rounded but feasible/valid version of the unconstrained solution $\hat{{\b{U}}}_\augp(k)$), or both~\cite{dorfling_2020}. Recall that the dependence of $\b{p}(\ell)$ on $\Delta \b{u}_\ph(\ell)$ has to be imposed to obtain a valid/feasible initial solution. The initial sphere radius $\rho_\ini^2(k)$ corresponds to the cost of the initial solution $\b{U}_\augp^\ini(k)$.

\vspace{-.1cm}
\section{Matrices for \ac{c2} and \ac{c3}}\label{app:matrices}
We first consider \ac{c2}. Let $\b{\Pi}_\augp\in \mathbb{R}^{6 N_\p \times 6 N_\p}$ be defined by \looseness=-1
{\begin{equation*}
    \b{\Pi}_\augp = \begin{bmatrix}
        \I_3 &&&&& \\
        \0_3 & \I_3 &&&& \\
        -\I_3 & \0_3 & \I_3 &&& \\
        & \ddots & \ddots & \ddots && \\
        && -\I_3 & \0_3 & \I_3 & \\
        &&& \0_3 & \0_3 & \I_3  
    \end{bmatrix},
\end{equation*}}
where $\I_n \in \mathbb{R}^{n\times n}$ denotes the identity matrix. Further, we define $\0_n \in \mathbb{R}^{n\times n}$ and $\0_{n\times m} \in \mathbb{R}^{n\times m}$ as zero matrices. Notice that $\b{\Pi}_\augp$ is full-rank. With this definition we obtain $\b{P}(k) = [\Delta \b{u}_\ph^\top(k)\ \b{p}^\top(k) \ldots \Delta \b{u}_\ph^\top(k+N_\p-1)\ \b{p}^\top(k+N_\p-1)]^\top$ from 
$$\b{P}(k) = \b{\Pi}_\augp\b{U}_\augp(k) - \b{E}_\augp\b{u}_\augp(k-1),$$ 
given that 
$$\b{E}_\augp = \begin{bmatrix}
    \I_3 & \0_3 & \cdots & \0_3\\
    \0_3 & \0_3 & \cdots & \0_3
\end{bmatrix}^\top \in \mathbb{R}^{6N_\p \times 6}.$$ 

For the \ac{c3}, we define $\b{\Pi}_\augs\in \mathbb{R}^{4 N_\p \times 4 N_\p}$ and $\b{L}_\augs \in \mathbb{R}^{N_\p \times 4 N_\p}$. The matrix $\b{L}_\augs$ is defined to extract the slack variable, i.e., $\b{S}(k+1) = \b{L}_\augs \b{U}_\augs(k)$. The matrix $\b{\Pi}_\augs$ is defined as
{\begin{equation*}
\b{\Pi}_\augs = \begin{bmatrix}
        \I_3 &&&&& \\
        \0_{1\times 3} & 0 &&&& \\
        -\I_3 & \0_{3\times 1} & \I_3 &&& \\
        & \ddots & \ddots & \ddots && \\
        && -\I_3 & \0_{3\times 1} & \I_3 & \\
        &&& 0 & \0_{1\times 3} & 0  
    \end{bmatrix},
\end{equation*}}
so that 
\begin{equation*}
    \Delta \b{U}_\augs(k) = \b{\Pi}_\augs(k)\b{U}_\augs(k) - \b{E}_\augs \b{u}_\augs(k-1)
\end{equation*}
hold,s given that 
\begin{equation*}
    \b{E}_\augs = \begin{bmatrix}
        \I_3 & \0_{3\times 1} & \0_{3\times 4(N_\p-1)} \\
        \0_{1\times 3} & 0 & \0_{1 \times 4(N_\p-1)}
    \end{bmatrix}^\top
    \in \mathbb{R}^{4 N_\p \times 4}.
\end{equation*}
Notice that $\lambda_\U \b{\Pi}_\augs^\top \b{\Pi}_\augs + \lambda_\sw \b{L}_\augs^2 \b{L}_\augs \succ 0$ is satisfied for all $\lambda_\U, \lambda_\sw > 0$.

Let $\bar{\b{Q}}_\augp = \I_{N_\p} \otimes \b{Q}_\augp \in \mathbb{R}^{3N_\p\times3N_\p}$ be the full-horizon output weight of \ac{c2}, where $\otimes$ denotes the Kronecker product. The output weight of \ac{c3} can similarly be defined as $\bar{\b{Q}}_\augs = \I_{N_\p} \otimes \b{Q}_\augs = \I_{2N_\p}$.\looseness=-1

For both  \ac{c2} and  \ac{c3}, the system dynamics are represented with the matrices $\b{A}_{(\cdot)} \in \mathbb{R}^{n\times n}$, $\b{B}_{(\cdot)} \in \mathbb{R}^{n\times m}$, $\b{D}_{(\cdot)} \in \mathbb{R}^{n\times 3}$, and $\b{C}_{(\cdot)} \in \mathbb{R}^{q\times n}$, where $(\cdot)\in\{\augp,\augs\}$, $q=3$, $n=4$, $m=6$ for \ac{c2}, and $q=2$, $n=2$, $m=4$ for \ac{c3}. The prediction horizon of each controller can be denoted as
$\b{Y}_{(\cdot)}(k+1) = [\b{y}_{(\cdot)}(k+1)^\top\ \ldots\ \b{y}_{(\cdot)}(k+N_\p)^\top]^\top,$
which results in the dynamics
\begin{equation*}
    \b{Y}_{(\cdot)}(k+1) = \b{\Gamma}_{(\cdot)} \b{x}_{(\cdot)}(k) + \b{\Upsilon}_{(\cdot)} \b{U}_{(\cdot)}(k) + \b{\Psi}_{(\cdot)} \b{V}_\g(k),
\end{equation*}
where the following definitions hold:
\begin{equation*}
    \b{\Gamma}_{(\cdot)} = \left[\begin{array}{c}
        % \I \\
        \b{C}_{(\cdot)}\b{A}_{(\cdot)} \\
        \vdots \\
        \b{C}_{(\cdot)}\b{A}_{(\cdot)}^{N_\p}
    \end{array}\right],
\end{equation*}    
    \begin{equation*}
\begin{split}
    \b{\Upsilon}_{(\cdot)} &= \left[\begin{array}{cccc}
        % \0 & & \ldots & \0 \\
        \b{C}_{(\cdot)} \b{B}_{(\cdot)} & \0_{q\times m} & \cdots & \0_{q\times m} \\
        \b{C}_{(\cdot)} \b{A}_{(\cdot)}\b{B}_{(\cdot)} & \b{C}_{(\cdot)} \b{B}_{(\cdot)} & & \vdots \\
        \vdots & & \ddots & \0_{q\times m} \\
        \b{C}_{(\cdot)}\b{A}_{(\cdot)}^{N_\p-1}\b{B}_{(\cdot)} & & \ldots & \b{C}_{(\cdot)} \b{B}_{(\cdot)}
    \end{array}\right],\\
    \b{\Psi}_{(\cdot)} &= \left[\begin{array}{cccc}
        % \0 & & \ldots & \0 \\
        \b{C}_{(\cdot)} \b{D}_{(\cdot)} & \0_{q\times 3} & \cdots & \0_{q\times 3} \\
        \b{C}_{(\cdot)} \b{A}_{(\cdot)}\b{D}_{(\cdot)} & \b{C}_{(\cdot)} \b{D}_{(\cdot)} & & \vdots \\
        \vdots & & \ddots & \0_{q\times 3} \\
        \b{C}_{(\cdot)} \b{A}_{(\cdot)}^{N_\p-1}\b{D}_{(\cdot)} & & \ldots & \b{C}_{(\cdot)} \b{D}_{(\cdot)}
    \end{array}\right].
    \end{split}
\end{equation*}
The matrices are of the sizes $\b{\Gamma}_{(\cdot)} \in \mathbb{R}^{q N_\p \times n}$, $\b{\Upsilon}_{(\cdot)} \in \mathbb{R}^{q N_\p \times m N_\p}$, and $\b{\Psi}_{(\cdot)} \in \mathbb{R}^{q N_\p \times 3 N_\p}$.
    
\end{document}

%% file: Figures/npc.tex
\begin{tikzpicture}[scale=1.5]%
% dpic version 2021.04.10 option -g for TikZ and PGF 1.01
\ifx\dpiclw\undefined\newdimen\dpiclw\fi
\global\def\dpicdraw{\draw[line width=\dpiclw]}
\global\def\dpicstop{;}
\dpiclw=0.8bp
\dpiclw=0.8bp
\footnotesize
\dpicdraw (0.401253,1.01)
 --(0.401253,0.883584)\dpicstop
\dpicdraw (0.401253,0.883584)
 --(0.44292,0.883584)
 --(0.401253,0.81695)
 --(0.359587,0.883584)
 --(0.401253,0.883584)\dpicstop
\dpicdraw (0.355535,0.811416)
 --(0.446972,0.811416)\dpicstop
\dpicdraw (0.38042,0.811416)
 --(0.338753,0.739247)\dpicstop
\dpicdraw (0.401253,0.811416)
 --(0.401253,0.685)\dpicstop
\dpicdraw (0.701253,0.685)
 --(0.701253,0.811416)\dpicstop
\global\let\dpicshdraw=\dpicdraw\global\def\dpicdraw{}
\global\def\dpicstop{--}
\dpicshdraw[fill=white!0!black]
\dpicdraw (0.701253,0.811416)
 --(0.659587,0.811416)
 --(0.701253,0.87805)
 --(0.74292,0.811416)
 --(0.701253,0.811416)\dpicstop
cycle; \global\let\dpicdraw=\dpicshdraw\global\def\dpicstop{;}
\dpicdraw (0.746972,0.883584)
 --(0.655535,0.883584)\dpicstop
\dpicdraw (0.701253,0.883584)
 --(0.701253,1.01)\dpicstop
\dpicdraw[line width=0.4bp](0.401253,1.01) circle (0.00109in)\dpicstop
\dpicdraw[line width=0.4bp](0.401253,0.685) circle (0.00109in)\dpicstop
\dpicdraw (0.401253,1.01)
 --(0.701253,1.01)\dpicstop
\dpicdraw[line width=0.4bp](0.701253,1.01) circle (0.00109in)\dpicstop
\dpicdraw (0.401253,0.685)
 --(0.701253,0.685)\dpicstop
\dpicdraw[line width=0.4bp](0.701253,0.685) circle (0.00109in)\dpicstop
\dpicdraw[fill=black](0.551253,1.01) circle (0.007874in)\dpicstop
\dpicdraw[fill=black](0.551253,0.685) circle (0.007874in)\dpicstop
\dpicdraw (0.401253,0.535)
 --(0.401253,0.408584)\dpicstop
\dpicdraw (0.401253,0.408584)
 --(0.44292,0.408584)
 --(0.401253,0.34195)
 --(0.359587,0.408584)
 --(0.401253,0.408584)\dpicstop
\dpicdraw (0.355535,0.336416)
 --(0.446972,0.336416)\dpicstop
\dpicdraw (0.38042,0.336416)
 --(0.338753,0.264247)\dpicstop
\dpicdraw (0.401253,0.336416)
 --(0.401253,0.21)\dpicstop
\dpicdraw (0.701253,0.21)
 --(0.701253,0.336416)\dpicstop
\global\let\dpicshdraw=\dpicdraw\global\def\dpicdraw{}
\global\def\dpicstop{--}
\dpicshdraw[fill=white!0!black]
\dpicdraw (0.701253,0.336416)
 --(0.659587,0.336416)
 --(0.701253,0.40305)
 --(0.74292,0.336416)
 --(0.701253,0.336416)\dpicstop
cycle; \global\let\dpicdraw=\dpicshdraw\global\def\dpicstop{;}
\dpicdraw (0.746972,0.408584)
 --(0.655535,0.408584)\dpicstop
\dpicdraw (0.701253,0.408584)
 --(0.701253,0.535)\dpicstop
\dpicdraw[line width=0.4bp](0.401253,0.535) circle (0.00109in)\dpicstop
\dpicdraw[line width=0.4bp](0.401253,0.21) circle (0.00109in)\dpicstop
\dpicdraw (0.401253,0.535)
 --(0.701253,0.535)\dpicstop
\dpicdraw[line width=0.4bp](0.701253,0.535) circle (0.00109in)\dpicstop
\dpicdraw (0.401253,0.21)
 --(0.701253,0.21)\dpicstop
\dpicdraw[line width=0.4bp](0.701253,0.21) circle (0.00109in)\dpicstop
\dpicdraw[fill=black](0.551253,0.535) circle (0.007874in)\dpicstop
\dpicdraw[fill=black](0.551253,0.21) circle (0.007874in)\dpicstop
\dpicdraw (0.401253,-0.21)
 --(0.401253,-0.336416)\dpicstop
\dpicdraw (0.401253,-0.336416)
 --(0.44292,-0.336416)
 --(0.401253,-0.40305)
 --(0.359587,-0.336416)
 --(0.401253,-0.336416)\dpicstop
\dpicdraw (0.355535,-0.408584)
 --(0.446972,-0.408584)\dpicstop
\dpicdraw (0.38042,-0.408584)
 --(0.338753,-0.480753)\dpicstop
\dpicdraw (0.401253,-0.408584)
 --(0.401253,-0.535)\dpicstop
\dpicdraw (0.701253,-0.535)
 --(0.701253,-0.408584)\dpicstop
\global\let\dpicshdraw=\dpicdraw\global\def\dpicdraw{}
\global\def\dpicstop{--}
\dpicshdraw[fill=white!0!black]
\dpicdraw (0.701253,-0.408584)
 --(0.659587,-0.408584)
 --(0.701253,-0.34195)
 --(0.74292,-0.408584)
 --(0.701253,-0.408584)\dpicstop
cycle; \global\let\dpicdraw=\dpicshdraw\global\def\dpicstop{;}
\dpicdraw (0.746972,-0.336416)
 --(0.655535,-0.336416)\dpicstop
\dpicdraw (0.701253,-0.336416)
 --(0.701253,-0.21)\dpicstop
\dpicdraw[line width=0.4bp](0.401253,-0.21) circle (0.00109in)\dpicstop
\dpicdraw[line width=0.4bp](0.401253,-0.535) circle (0.00109in)\dpicstop
\dpicdraw (0.401253,-0.21)
 --(0.701253,-0.21)\dpicstop
\dpicdraw[line width=0.4bp](0.701253,-0.21) circle (0.00109in)\dpicstop
\dpicdraw (0.401253,-0.535)
 --(0.701253,-0.535)\dpicstop
\dpicdraw[line width=0.4bp](0.701253,-0.535) circle (0.00109in)\dpicstop
\dpicdraw[fill=black](0.551253,-0.21) circle (0.007874in)\dpicstop
\dpicdraw[fill=black](0.551253,-0.535) circle (0.007874in)\dpicstop
\dpicdraw (0.401253,-0.685)
 --(0.401253,-0.811416)\dpicstop
\dpicdraw (0.401253,-0.811416)
 --(0.44292,-0.811416)
 --(0.401253,-0.87805)
 --(0.359587,-0.811416)
 --(0.401253,-0.811416)\dpicstop
\dpicdraw (0.355535,-0.883584)
 --(0.446972,-0.883584)\dpicstop
\dpicdraw (0.38042,-0.883584)
 --(0.338753,-0.955753)\dpicstop
\dpicdraw (0.401253,-0.883584)
 --(0.401253,-1.01)\dpicstop
\dpicdraw (0.701253,-1.01)
 --(0.701253,-0.883584)\dpicstop
\global\let\dpicshdraw=\dpicdraw\global\def\dpicdraw{}
\global\def\dpicstop{--}
\dpicshdraw[fill=white!0!black]
\dpicdraw (0.701253,-0.883584)
 --(0.659587,-0.883584)
 --(0.701253,-0.81695)
 --(0.74292,-0.883584)
 --(0.701253,-0.883584)\dpicstop
cycle; \global\let\dpicdraw=\dpicshdraw\global\def\dpicstop{;}
\dpicdraw (0.746972,-0.811416)
 --(0.655535,-0.811416)\dpicstop
\dpicdraw (0.701253,-0.811416)
 --(0.701253,-0.685)\dpicstop
\dpicdraw[line width=0.4bp](0.401253,-0.685) circle (0.00109in)\dpicstop
\dpicdraw[line width=0.4bp](0.401253,-1.01) circle (0.00109in)\dpicstop
\dpicdraw (0.401253,-0.685)
 --(0.701253,-0.685)\dpicstop
\dpicdraw[line width=0.4bp](0.701253,-0.685) circle (0.00109in)\dpicstop
\dpicdraw (0.401253,-1.01)
 --(0.701253,-1.01)\dpicstop
\dpicdraw[line width=0.4bp](0.701253,-1.01) circle (0.00109in)\dpicstop
\dpicdraw[fill=black](0.551253,-0.685) circle (0.007874in)\dpicstop
\dpicdraw[fill=black](0.551253,-1.01) circle (0.007874in)\dpicstop
\dpicdraw (0.551253,0.685)
 --(0.551253,0.535)\dpicstop
\dpicdraw (0.551253,0.21)
 --(0.551253,-0.21)\dpicstop
\dpicdraw (0.551253,-0.535)
 --(0.551253,-0.685)\dpicstop
\dpicdraw[fill=black](0.551253,0.61) circle (0.007874in)\dpicstop
\dpicdraw[fill=black](0.551253,-0.61) circle (0.007874in)\dpicstop
\dpicdraw (0.551253,0.61)
 --(0.051253,0.61)\dpicstop
\dpicdraw[line width=0.4bp](0.051253,0.61) circle (0.00109in)\dpicstop
\dpicdraw (0.051253,-0.084)
 --(0.051253,0.226916)\dpicstop
\global\let\dpicshdraw=\dpicdraw\global\def\dpicdraw{}
\global\def\dpicstop{--}
\dpicshdraw[fill=white!0!black]
\dpicdraw (0.051253,0.226916)
 --(0.009587,0.226916)
 --(0.051253,0.29355)
 --(0.09292,0.226916)
 --(0.051253,0.226916)\dpicstop
cycle; \global\let\dpicdraw=\dpicshdraw\global\def\dpicstop{;}
\dpicdraw (0.096972,0.299084)
 --(0.005535,0.299084)\dpicstop
\dpicdraw (0.051253,0.299084)
 --(0.051253,0.61)\dpicstop
\dpicdraw[fill=black](0.051253,-0.084) circle (0.007874in)\dpicstop
\dpicdraw (0.051253,-0.61)
 --(0.051253,-0.383084)\dpicstop
\global\let\dpicshdraw=\dpicdraw\global\def\dpicdraw{}
\global\def\dpicstop{--}
\dpicshdraw[fill=white!0!black]
\dpicdraw (0.051253,-0.383084)
 --(0.009587,-0.383084)
 --(0.051253,-0.31645)
 --(0.09292,-0.383084)
 --(0.051253,-0.383084)\dpicstop
cycle; \global\let\dpicdraw=\dpicshdraw\global\def\dpicstop{;}
\dpicdraw (0.096972,-0.310916)
 --(0.005535,-0.310916)\dpicstop
\dpicdraw (0.051253,-0.310916)
 --(0.051253,-0.084)\dpicstop
\dpicdraw[line width=0.4bp](0.051253,-0.61) circle (0.00109in)\dpicstop
\dpicdraw (0.051253,-0.61)
 --(0.551253,-0.61)\dpicstop
\dpicdraw (1.451253,1.01)
 --(1.451253,0.883584)\dpicstop
\dpicdraw (1.451253,0.883584)
 --(1.49292,0.883584)
 --(1.451253,0.81695)
 --(1.409587,0.883584)
 --(1.451253,0.883584)\dpicstop
\dpicdraw (1.405535,0.811416)
 --(1.496972,0.811416)\dpicstop
\dpicdraw (1.43042,0.811416)
 --(1.388753,0.739247)\dpicstop
\dpicdraw (1.451253,0.811416)
 --(1.451253,0.685)\dpicstop
\dpicdraw (1.751253,0.685)
 --(1.751253,0.811416)\dpicstop
\global\let\dpicshdraw=\dpicdraw\global\def\dpicdraw{}
\global\def\dpicstop{--}
\dpicshdraw[fill=white!0!black]
\dpicdraw (1.751253,0.811416)
 --(1.709587,0.811416)
 --(1.751253,0.87805)
 --(1.79292,0.811416)
 --(1.751253,0.811416)\dpicstop
cycle; \global\let\dpicdraw=\dpicshdraw\global\def\dpicstop{;}
\dpicdraw (1.796972,0.883584)
 --(1.705535,0.883584)\dpicstop
\dpicdraw (1.751253,0.883584)
 --(1.751253,1.01)\dpicstop
\dpicdraw[line width=0.4bp](1.451253,1.01) circle (0.00109in)\dpicstop
\dpicdraw[line width=0.4bp](1.451253,0.685) circle (0.00109in)\dpicstop
\dpicdraw (1.451253,1.01)
 --(1.751253,1.01)\dpicstop
\dpicdraw[line width=0.4bp](1.751253,1.01) circle (0.00109in)\dpicstop
\dpicdraw (1.451253,0.685)
 --(1.751253,0.685)\dpicstop
\dpicdraw[line width=0.4bp](1.751253,0.685) circle (0.00109in)\dpicstop
\dpicdraw[fill=black](1.601253,1.01) circle (0.007874in)\dpicstop
\dpicdraw[fill=black](1.601253,0.685) circle (0.007874in)\dpicstop
\dpicdraw (1.451253,0.535)
 --(1.451253,0.408584)\dpicstop
\dpicdraw (1.451253,0.408584)
 --(1.49292,0.408584)
 --(1.451253,0.34195)
 --(1.409587,0.408584)
 --(1.451253,0.408584)\dpicstop
\dpicdraw (1.405535,0.336416)
 --(1.496972,0.336416)\dpicstop
\dpicdraw (1.43042,0.336416)
 --(1.388753,0.264247)\dpicstop
\dpicdraw (1.451253,0.336416)
 --(1.451253,0.21)\dpicstop
\dpicdraw (1.751253,0.21)
 --(1.751253,0.336416)\dpicstop
\global\let\dpicshdraw=\dpicdraw\global\def\dpicdraw{}
\global\def\dpicstop{--}
\dpicshdraw[fill=white!0!black]
\dpicdraw (1.751253,0.336416)
 --(1.709587,0.336416)
 --(1.751253,0.40305)
 --(1.79292,0.336416)
 --(1.751253,0.336416)\dpicstop
cycle; \global\let\dpicdraw=\dpicshdraw\global\def\dpicstop{;}
\dpicdraw (1.796972,0.408584)
 --(1.705535,0.408584)\dpicstop
\dpicdraw (1.751253,0.408584)
 --(1.751253,0.535)\dpicstop
\dpicdraw[line width=0.4bp](1.451253,0.535) circle (0.00109in)\dpicstop
\dpicdraw[line width=0.4bp](1.451253,0.21) circle (0.00109in)\dpicstop
\dpicdraw (1.451253,0.535)
 --(1.751253,0.535)\dpicstop
\dpicdraw[line width=0.4bp](1.751253,0.535) circle (0.00109in)\dpicstop
\dpicdraw (1.451253,0.21)
 --(1.751253,0.21)\dpicstop
\dpicdraw[line width=0.4bp](1.751253,0.21) circle (0.00109in)\dpicstop
\dpicdraw[fill=black](1.601253,0.535) circle (0.007874in)\dpicstop
\dpicdraw[fill=black](1.601253,0.21) circle (0.007874in)\dpicstop
\dpicdraw (1.451253,-0.21)
 --(1.451253,-0.336416)\dpicstop
\dpicdraw (1.451253,-0.336416)
 --(1.49292,-0.336416)
 --(1.451253,-0.40305)
 --(1.409587,-0.336416)
 --(1.451253,-0.336416)\dpicstop
\dpicdraw (1.405535,-0.408584)
 --(1.496972,-0.408584)\dpicstop
\dpicdraw (1.43042,-0.408584)
 --(1.388753,-0.480753)\dpicstop
\dpicdraw (1.451253,-0.408584)
 --(1.451253,-0.535)\dpicstop
\dpicdraw (1.751253,-0.535)
 --(1.751253,-0.408584)\dpicstop
\global\let\dpicshdraw=\dpicdraw\global\def\dpicdraw{}
\global\def\dpicstop{--}
\dpicshdraw[fill=white!0!black]
\dpicdraw (1.751253,-0.408584)
 --(1.709587,-0.408584)
 --(1.751253,-0.34195)
 --(1.79292,-0.408584)
 --(1.751253,-0.408584)\dpicstop
cycle; \global\let\dpicdraw=\dpicshdraw\global\def\dpicstop{;}
\dpicdraw (1.796972,-0.336416)
 --(1.705535,-0.336416)\dpicstop
\dpicdraw (1.751253,-0.336416)
 --(1.751253,-0.21)\dpicstop
\dpicdraw[line width=0.4bp](1.451253,-0.21) circle (0.00109in)\dpicstop
\dpicdraw[line width=0.4bp](1.451253,-0.535) circle (0.00109in)\dpicstop
\dpicdraw (1.451253,-0.21)
 --(1.751253,-0.21)\dpicstop
\dpicdraw[line width=0.4bp](1.751253,-0.21) circle (0.00109in)\dpicstop
\dpicdraw (1.451253,-0.535)
 --(1.751253,-0.535)\dpicstop
\dpicdraw[line width=0.4bp](1.751253,-0.535) circle (0.00109in)\dpicstop
\dpicdraw[fill=black](1.601253,-0.21) circle (0.007874in)\dpicstop
\dpicdraw[fill=black](1.601253,-0.535) circle (0.007874in)\dpicstop
\dpicdraw (1.451253,-0.685)
 --(1.451253,-0.811416)\dpicstop
\dpicdraw (1.451253,-0.811416)
 --(1.49292,-0.811416)
 --(1.451253,-0.87805)
 --(1.409587,-0.811416)
 --(1.451253,-0.811416)\dpicstop
\dpicdraw (1.405535,-0.883584)
 --(1.496972,-0.883584)\dpicstop
\dpicdraw (1.43042,-0.883584)
 --(1.388753,-0.955753)\dpicstop
\dpicdraw (1.451253,-0.883584)
 --(1.451253,-1.01)\dpicstop
\dpicdraw (1.751253,-1.01)
 --(1.751253,-0.883584)\dpicstop
\global\let\dpicshdraw=\dpicdraw\global\def\dpicdraw{}
\global\def\dpicstop{--}
\dpicshdraw[fill=white!0!black]
\dpicdraw (1.751253,-0.883584)
 --(1.709587,-0.883584)
 --(1.751253,-0.81695)
 --(1.79292,-0.883584)
 --(1.751253,-0.883584)\dpicstop
cycle; \global\let\dpicdraw=\dpicshdraw\global\def\dpicstop{;}
\dpicdraw (1.796972,-0.811416)
 --(1.705535,-0.811416)\dpicstop
\dpicdraw (1.751253,-0.811416)
 --(1.751253,-0.685)\dpicstop
\dpicdraw[line width=0.4bp](1.451253,-0.685) circle (0.00109in)\dpicstop
\dpicdraw[line width=0.4bp](1.451253,-1.01) circle (0.00109in)\dpicstop
\dpicdraw (1.451253,-0.685)
 --(1.751253,-0.685)\dpicstop
\dpicdraw[line width=0.4bp](1.751253,-0.685) circle (0.00109in)\dpicstop
\dpicdraw (1.451253,-1.01)
 --(1.751253,-1.01)\dpicstop
\dpicdraw[line width=0.4bp](1.751253,-1.01) circle (0.00109in)\dpicstop
\dpicdraw[fill=black](1.601253,-0.685) circle (0.007874in)\dpicstop
\dpicdraw[fill=black](1.601253,-1.01) circle (0.007874in)\dpicstop
\dpicdraw (1.601253,0.685)
 --(1.601253,0.535)\dpicstop
\dpicdraw (1.601253,0.21)
 --(1.601253,-0.21)\dpicstop
\dpicdraw (1.601253,-0.535)
 --(1.601253,-0.685)\dpicstop
\dpicdraw[fill=black](1.601253,0.61) circle (0.007874in)\dpicstop
\dpicdraw[fill=black](1.601253,-0.61) circle (0.007874in)\dpicstop
\dpicdraw (1.601253,0.61)
 --(1.101253,0.61)\dpicstop
\dpicdraw[line width=0.4bp](1.101253,0.61) circle (0.00109in)\dpicstop
\dpicdraw (1.101253,-0.084)
 --(1.101253,0.226916)\dpicstop
\global\let\dpicshdraw=\dpicdraw\global\def\dpicdraw{}
\global\def\dpicstop{--}
\dpicshdraw[fill=white!0!black]
\dpicdraw (1.101253,0.226916)
 --(1.059587,0.226916)
 --(1.101253,0.29355)
 --(1.14292,0.226916)
 --(1.101253,0.226916)\dpicstop
cycle; \global\let\dpicdraw=\dpicshdraw\global\def\dpicstop{;}
\dpicdraw (1.146972,0.299084)
 --(1.055535,0.299084)\dpicstop
\dpicdraw (1.101253,0.299084)
 --(1.101253,0.61)\dpicstop
\dpicdraw[fill=black](1.101253,-0.084) circle (0.007874in)\dpicstop
\dpicdraw (1.101253,-0.61)
 --(1.101253,-0.383084)\dpicstop
\global\let\dpicshdraw=\dpicdraw\global\def\dpicdraw{}
\global\def\dpicstop{--}
\dpicshdraw[fill=white!0!black]
\dpicdraw (1.101253,-0.383084)
 --(1.059587,-0.383084)
 --(1.101253,-0.31645)
 --(1.14292,-0.383084)
 --(1.101253,-0.383084)\dpicstop
cycle; \global\let\dpicdraw=\dpicshdraw\global\def\dpicstop{;}
\dpicdraw (1.146972,-0.310916)
 --(1.055535,-0.310916)\dpicstop
\dpicdraw (1.101253,-0.310916)
 --(1.101253,-0.084)\dpicstop
\dpicdraw[line width=0.4bp](1.101253,-0.61) circle (0.00109in)\dpicstop
\dpicdraw (1.101253,-0.61)
 --(1.601253,-0.61)\dpicstop
\dpicdraw (2.501253,1.01)
 --(2.501253,0.883584)\dpicstop
\dpicdraw (2.501253,0.883584)
 --(2.54292,0.883584)
 --(2.501253,0.81695)
 --(2.459587,0.883584)
 --(2.501253,0.883584)\dpicstop
\dpicdraw (2.455535,0.811416)
 --(2.546972,0.811416)\dpicstop
\dpicdraw (2.48042,0.811416)
 --(2.438753,0.739247)\dpicstop
\dpicdraw (2.501253,0.811416)
 --(2.501253,0.685)\dpicstop
\dpicdraw (2.801253,0.685)
 --(2.801253,0.811416)\dpicstop
\global\let\dpicshdraw=\dpicdraw\global\def\dpicdraw{}
\global\def\dpicstop{--}
\dpicshdraw[fill=white!0!black]
\dpicdraw (2.801253,0.811416)
 --(2.759587,0.811416)
 --(2.801253,0.87805)
 --(2.84292,0.811416)
 --(2.801253,0.811416)\dpicstop
cycle; \global\let\dpicdraw=\dpicshdraw\global\def\dpicstop{;}
\dpicdraw (2.846972,0.883584)
 --(2.755535,0.883584)\dpicstop
\dpicdraw (2.801253,0.883584)
 --(2.801253,1.01)\dpicstop
\dpicdraw[line width=0.4bp](2.501253,1.01) circle (0.00109in)\dpicstop
\dpicdraw[line width=0.4bp](2.501253,0.685) circle (0.00109in)\dpicstop
\dpicdraw (2.501253,1.01)
 --(2.801253,1.01)\dpicstop
\dpicdraw[line width=0.4bp](2.801253,1.01) circle (0.00109in)\dpicstop
\dpicdraw (2.501253,0.685)
 --(2.801253,0.685)\dpicstop
\dpicdraw[line width=0.4bp](2.801253,0.685) circle (0.00109in)\dpicstop
\dpicdraw[fill=black](2.651253,1.01) circle (0.007874in)\dpicstop
\dpicdraw[fill=black](2.651253,0.685) circle (0.007874in)\dpicstop
\dpicdraw (2.501253,0.535)
 --(2.501253,0.408584)\dpicstop
\dpicdraw (2.501253,0.408584)
 --(2.54292,0.408584)
 --(2.501253,0.34195)
 --(2.459587,0.408584)
 --(2.501253,0.408584)\dpicstop
\dpicdraw (2.455535,0.336416)
 --(2.546972,0.336416)\dpicstop
\dpicdraw (2.48042,0.336416)
 --(2.438753,0.264247)\dpicstop
\dpicdraw (2.501253,0.336416)
 --(2.501253,0.21)\dpicstop
\dpicdraw (2.801253,0.21)
 --(2.801253,0.336416)\dpicstop
\global\let\dpicshdraw=\dpicdraw\global\def\dpicdraw{}
\global\def\dpicstop{--}
\dpicshdraw[fill=white!0!black]
\dpicdraw (2.801253,0.336416)
 --(2.759587,0.336416)
 --(2.801253,0.40305)
 --(2.84292,0.336416)
 --(2.801253,0.336416)\dpicstop
cycle; \global\let\dpicdraw=\dpicshdraw\global\def\dpicstop{;}
\dpicdraw (2.846972,0.408584)
 --(2.755535,0.408584)\dpicstop
\dpicdraw (2.801253,0.408584)
 --(2.801253,0.535)\dpicstop
\dpicdraw[line width=0.4bp](2.501253,0.535) circle (0.00109in)\dpicstop
\dpicdraw[line width=0.4bp](2.501253,0.21) circle (0.00109in)\dpicstop
\dpicdraw (2.501253,0.535)
 --(2.801253,0.535)\dpicstop
\dpicdraw[line width=0.4bp](2.801253,0.535) circle (0.00109in)\dpicstop
\dpicdraw (2.501253,0.21)
 --(2.801253,0.21)\dpicstop
\dpicdraw[line width=0.4bp](2.801253,0.21) circle (0.00109in)\dpicstop
\dpicdraw[fill=black](2.651253,0.535) circle (0.007874in)\dpicstop
\dpicdraw[fill=black](2.651253,0.21) circle (0.007874in)\dpicstop
\dpicdraw (2.501253,-0.21)
 --(2.501253,-0.336416)\dpicstop
\dpicdraw (2.501253,-0.336416)
 --(2.54292,-0.336416)
 --(2.501253,-0.40305)
 --(2.459587,-0.336416)
 --(2.501253,-0.336416)\dpicstop
\dpicdraw (2.455535,-0.408584)
 --(2.546972,-0.408584)\dpicstop
\dpicdraw (2.48042,-0.408584)
 --(2.438753,-0.480753)\dpicstop
\dpicdraw (2.501253,-0.408584)
 --(2.501253,-0.535)\dpicstop
\dpicdraw (2.801253,-0.535)
 --(2.801253,-0.408584)\dpicstop
\global\let\dpicshdraw=\dpicdraw\global\def\dpicdraw{}
\global\def\dpicstop{--}
\dpicshdraw[fill=white!0!black]
\dpicdraw (2.801253,-0.408584)
 --(2.759587,-0.408584)
 --(2.801253,-0.34195)
 --(2.84292,-0.408584)
 --(2.801253,-0.408584)\dpicstop
cycle; \global\let\dpicdraw=\dpicshdraw\global\def\dpicstop{;}
\dpicdraw (2.846972,-0.336416)
 --(2.755535,-0.336416)\dpicstop
\dpicdraw (2.801253,-0.336416)
 --(2.801253,-0.21)\dpicstop
\dpicdraw[line width=0.4bp](2.501253,-0.21) circle (0.00109in)\dpicstop
\dpicdraw[line width=0.4bp](2.501253,-0.535) circle (0.00109in)\dpicstop
\dpicdraw (2.501253,-0.21)
 --(2.801253,-0.21)\dpicstop
\dpicdraw[line width=0.4bp](2.801253,-0.21) circle (0.00109in)\dpicstop
\dpicdraw (2.501253,-0.535)
 --(2.801253,-0.535)\dpicstop
\dpicdraw[line width=0.4bp](2.801253,-0.535) circle (0.00109in)\dpicstop
\dpicdraw[fill=black](2.651253,-0.21) circle (0.007874in)\dpicstop
\dpicdraw[fill=black](2.651253,-0.535) circle (0.007874in)\dpicstop
\dpicdraw (2.501253,-0.685)
 --(2.501253,-0.811416)\dpicstop
\dpicdraw (2.501253,-0.811416)
 --(2.54292,-0.811416)
 --(2.501253,-0.87805)
 --(2.459587,-0.811416)
 --(2.501253,-0.811416)\dpicstop
\dpicdraw (2.455535,-0.883584)
 --(2.546972,-0.883584)\dpicstop
\dpicdraw (2.48042,-0.883584)
 --(2.438753,-0.955753)\dpicstop
\dpicdraw (2.501253,-0.883584)
 --(2.501253,-1.01)\dpicstop
\dpicdraw (2.801253,-1.01)
 --(2.801253,-0.883584)\dpicstop
\global\let\dpicshdraw=\dpicdraw\global\def\dpicdraw{}
\global\def\dpicstop{--}
\dpicshdraw[fill=white!0!black]
\dpicdraw (2.801253,-0.883584)
 --(2.759587,-0.883584)
 --(2.801253,-0.81695)
 --(2.84292,-0.883584)
 --(2.801253,-0.883584)\dpicstop
cycle; \global\let\dpicdraw=\dpicshdraw\global\def\dpicstop{;}
\dpicdraw (2.846972,-0.811416)
 --(2.755535,-0.811416)\dpicstop
\dpicdraw (2.801253,-0.811416)
 --(2.801253,-0.685)\dpicstop
\dpicdraw[line width=0.4bp](2.501253,-0.685) circle (0.00109in)\dpicstop
\dpicdraw[line width=0.4bp](2.501253,-1.01) circle (0.00109in)\dpicstop
\dpicdraw (2.501253,-0.685)
 --(2.801253,-0.685)\dpicstop
\dpicdraw[line width=0.4bp](2.801253,-0.685) circle (0.00109in)\dpicstop
\dpicdraw (2.501253,-1.01)
 --(2.801253,-1.01)\dpicstop
\dpicdraw[line width=0.4bp](2.801253,-1.01) circle (0.00109in)\dpicstop
\dpicdraw[fill=black](2.651253,-0.685) circle (0.007874in)\dpicstop
\dpicdraw[fill=black](2.651253,-1.01) circle (0.007874in)\dpicstop
\dpicdraw (2.651253,0.685)
 --(2.651253,0.535)\dpicstop
\dpicdraw (2.651253,0.21)
 --(2.651253,-0.21)\dpicstop
\dpicdraw (2.651253,-0.535)
 --(2.651253,-0.685)\dpicstop
\dpicdraw[fill=black](2.651253,0.61) circle (0.007874in)\dpicstop
\dpicdraw[fill=black](2.651253,-0.61) circle (0.007874in)\dpicstop
\dpicdraw (2.651253,0.61)
 --(2.151253,0.61)\dpicstop
\dpicdraw[line width=0.4bp](2.151253,0.61) circle (0.00109in)\dpicstop
\dpicdraw (2.151253,-0.084)
 --(2.151253,0.226916)\dpicstop
\global\let\dpicshdraw=\dpicdraw\global\def\dpicdraw{}
\global\def\dpicstop{--}
\dpicshdraw[fill=white!0!black]
\dpicdraw (2.151253,0.226916)
 --(2.109587,0.226916)
 --(2.151253,0.29355)
 --(2.19292,0.226916)
 --(2.151253,0.226916)\dpicstop
cycle; \global\let\dpicdraw=\dpicshdraw\global\def\dpicstop{;}
\dpicdraw (2.196972,0.299084)
 --(2.105535,0.299084)\dpicstop
\dpicdraw (2.151253,0.299084)
 --(2.151253,0.61)\dpicstop
\dpicdraw[fill=black](2.151253,-0.084) circle (0.007874in)\dpicstop
\dpicdraw (2.151253,-0.61)
 --(2.151253,-0.383084)\dpicstop
\global\let\dpicshdraw=\dpicdraw\global\def\dpicdraw{}
\global\def\dpicstop{--}
\dpicshdraw[fill=white!0!black]
\dpicdraw (2.151253,-0.383084)
 --(2.109587,-0.383084)
 --(2.151253,-0.31645)
 --(2.19292,-0.383084)
 --(2.151253,-0.383084)\dpicstop
cycle; \global\let\dpicdraw=\dpicshdraw\global\def\dpicstop{;}
\dpicdraw (2.196972,-0.310916)
 --(2.105535,-0.310916)\dpicstop
\dpicdraw (2.151253,-0.310916)
 --(2.151253,-0.084)\dpicstop
\dpicdraw[line width=0.4bp](2.151253,-0.61) circle (0.00109in)\dpicstop
\dpicdraw (2.151253,-0.61)
 --(2.651253,-0.61)\dpicstop
\dpicdraw[fill=black](0.551253,0.105) circle (0.007874in)\dpicstop
\dpicdraw (0.551253,0.105)
 --(1.065142,0.105)\dpicstop
\dpicdraw (1.059587,0.105)
 ..controls (1.059587,0.160556) and (1.14292,0.160556)
 ..(1.14292,0.105)\dpicstop
\dpicdraw (1.137364,0.105)
 --(1.565142,0.105)\dpicstop
\dpicdraw (1.559587,0.105)
 ..controls (1.559587,0.160556) and (1.64292,0.160556)
 ..(1.64292,0.105)\dpicstop
\dpicdraw (1.637364,0.105)
 --(2.115142,0.105)\dpicstop
\dpicdraw (2.109587,0.105)
 ..controls (2.109587,0.160556) and (2.19292,0.160556)
 ..(2.19292,0.105)\dpicstop
\dpicdraw (2.187364,0.105)
 --(2.615142,0.105)\dpicstop
\dpicdraw (2.609587,0.105)
 ..controls (2.609587,0.160556) and (2.69292,0.160556)
 ..(2.69292,0.105)\dpicstop
\dpicdraw (2.687364,0.105)
 --(2.951253,0.105)\dpicstop
\dpicdraw[fill=black](1.601253,0) circle (0.007874in)\dpicstop
\dpicdraw (1.601253,0)
 --(2.115142,0)\dpicstop
\dpicdraw (2.109587,0)
 ..controls (2.109587,0.023012) and (2.128241,0.041667)
 ..(2.151253,0.041667)
 ..controls (2.174265,0.041667) and (2.19292,0.023012)
 ..(2.19292,-0)\dpicstop
\dpicdraw (2.187364,-0)
 --(2.615142,0)\dpicstop
\dpicdraw (2.609587,0)
 ..controls (2.609587,0.023012) and (2.628241,0.041667)
 ..(2.651253,0.041667)
 ..controls (2.674265,0.041667) and (2.69292,0.023012)
 ..(2.69292,-0)\dpicstop
\dpicdraw (2.687364,-0)
 --(2.951253,0)\dpicstop
\dpicdraw[fill=black](2.651253,-0.105) circle (0.007874in)\dpicstop
\dpicdraw (2.651253,-0.105)
 --(2.951253,-0.105)\dpicstop
\dpicdraw (2.651253,1.01)
 --(2.651253,1.115)\dpicstop
\dpicdraw[line width=0.4bp](2.651253,1.115) circle (0.00109in)\dpicstop
\dpicdraw (2.651253,1.115)
 --(-0.348747,1.115)\dpicstop
\dpicdraw[line width=0.4bp](-0.348747,1.115) circle (0.00109in)\dpicstop
\dpicdraw (2.651253,-1.01)
 --(2.651253,-1.115)\dpicstop
\dpicdraw[line width=0.4bp](2.651253,-1.115) circle (0.00109in)\dpicstop
\dpicdraw (2.651253,-1.115)
 --(-0.348747,-1.115)\dpicstop
\dpicdraw[line width=0.4bp](-0.348747,-1.115) circle (0.00109in)\dpicstop
\dpicdraw (-0.348747,-0.084)
 --(-0.348747,0.3905)\dpicstop
\dpicdraw (-0.348747,0.5155) circle (0.049213in)\dpicstop
\draw (-0.348747,0.453) node{$_-$};
\draw (-0.348747,0.578) node{$_+$};
\dpicdraw (-0.348747,0.6405)
 --(-0.348747,1.115)\dpicstop
\draw (-0.348747,0.6405) node[above left=-2bp]{$ +$};
\draw (-0.473747,0.5155) node[left=-2bp]{$ \frac{V_{\mathrm{dc}}}{2}$};
\draw (-0.348747,0.3905) node[below left=-2bp]{$ -$};
\dpicdraw (-0.348747,-1.115)
 --(-0.348747,-0.7245)\dpicstop
\dpicdraw (-0.348747,-0.5995) circle (0.049213in)\dpicstop
\draw (-0.348747,-0.662) node{$_-$};
\draw (-0.348747,-0.537) node{$_+$};
\dpicdraw (-0.348747,-0.4745)
 --(-0.348747,-0.084)\dpicstop
\draw (-0.348747,-0.4745) node[above left=-2bp]{$ +$};
\draw (-0.473747,-0.5995) node[left=-2bp]{$ \frac{V_{\mathrm{dc}}}{2}$};
\draw (-0.348747,-0.7245) node[below left=-2bp]{$ -$};
\dpicdraw[fill=black](-0.348747,-0.084) circle (0.007874in)\dpicstop
\dpicdraw (-0.348747,-0.084)
 --(0.515142,-0.084)\dpicstop
\dpicdraw (0.509587,-0.084)
 ..controls (0.509587,-0.028444) and (0.59292,-0.028444)
 ..(0.59292,-0.084)\dpicstop
\dpicdraw (0.587364,-0.084)
 --(1.565142,-0.084)\dpicstop
\dpicdraw (1.559587,-0.084)
 ..controls (1.559587,-0.028444) and (1.64292,-0.028444)
 ..(1.64292,-0.084)\dpicstop
\dpicdraw (1.637364,-0.084)
 --(2.151253,-0.084)\dpicstop
\dpicdraw (0.551253,1.01)
 --(0.551253,1.115)\dpicstop
\dpicdraw[fill=black](0.551253,1.115) circle (0.007874in)\dpicstop
\dpicdraw (1.601253,1.01)
 --(1.601253,1.115)\dpicstop
\dpicdraw[fill=black](1.601253,1.115) circle (0.007874in)\dpicstop
\dpicdraw (0.551253,-1.01)
 --(0.551253,-1.115)\dpicstop
\dpicdraw[fill=black](0.551253,-1.115) circle (0.007874in)\dpicstop
\dpicdraw (1.601253,-1.01)
 --(1.601253,-1.115)\dpicstop
\dpicdraw[fill=black](1.601253,-1.115) circle (0.007874in)\dpicstop
\dpicdraw[line width=0.4bp](2.951253,0.105) circle (0.00109in)\dpicstop
\dpicdraw (2.951253,0.105)
 --(2.951253,0.405)\dpicstop
\dpicdraw[line width=0.4bp](2.951253,0.405) circle (0.00109in)\dpicstop
\dpicdraw (2.951253,0.405)
 --(3.151253,0.405)\dpicstop
\draw (3.151253,0.405) node[above=-2bp]{$v_{\C,\A}$};
\dpicdraw (3.151253,0.405)
 --(3.276253,0.405)\dpicstop
\dpicdraw (3.276253,0.405)
 --(3.276253,0.399444)\dpicstop
\dpicdraw (3.276253,0.405)
 ..controls (3.276253,0.422259) and (3.290244,0.43625)
 ..(3.307503,0.43625)
 ..controls (3.324762,0.43625) and (3.338753,0.422259)
 ..(3.338753,0.405)\dpicstop
\dpicdraw (3.338753,0.405)
 --(3.338753,0.399444)\dpicstop
\dpicdraw (3.338753,0.405)
 ..controls (3.338753,0.422259) and (3.352744,0.43625)
 ..(3.370003,0.43625)
 ..controls (3.387262,0.43625) and (3.401253,0.422259)
 ..(3.401253,0.405)\dpicstop
\dpicdraw (3.401253,0.405)
 --(3.401253,0.399444)\dpicstop
\dpicdraw (3.401253,0.405)
 ..controls (3.401253,0.422259) and (3.415244,0.43625)
 ..(3.432503,0.43625)
 ..controls (3.449762,0.43625) and (3.463753,0.422259)
 ..(3.463753,0.405)\dpicstop
\dpicdraw (3.463753,0.405)
 --(3.463753,0.399444)\dpicstop
\dpicdraw (3.463753,0.405)
 ..controls (3.463753,0.422259) and (3.477744,0.43625)
 ..(3.495003,0.43625)
 ..controls (3.512262,0.43625) and (3.526253,0.422259)
 ..(3.526253,0.405)\dpicstop
\dpicdraw (3.526253,0.405)
 --(3.526253,0.399444)\dpicstop
\dpicdraw (3.526253,0.405)
 --(3.651253,0.405)\dpicstop
\filldraw (3.576253,0.38625)
 --(3.651253,0.405)
 --(3.576253,0.42375) --cycle\dpicstop
\dpicdraw (3.628347,0.405)
 --(3.576253,0.405)\dpicstop
\draw (3.6023,0.405) node[above=-2bp]{$ i_a$};
\draw (3.401253,0.43625) node[above=-2bp]{$ L$};
\dpicdraw[fill=white](3.151253,0.405) circle (0.007874in)\dpicstop
\dpicdraw (3.651253,0.405)
 --(3.776253,0.405)
 --(3.797087,0.446667)
 --(3.838753,0.363333)
 --(3.88042,0.446667)
 --(3.922087,0.363333)
 --(3.963753,0.446667)
 --(4.00542,0.363333)
 --(4.026253,0.405)
 --(4.151253,0.405)\dpicstop
\draw (3.901253,0.446667) node[above=-2bp]{$ R$};
\dpicdraw (4.151253,0.405)
 --(4.276253,0.405)\dpicstop
\dpicdraw (4.401253,0.405) circle (0.049213in)\dpicstop
\dpicdraw (4.401253,0.405)
 ..controls (4.401253,0.428012) and (4.382598,0.446667)
 ..(4.359587,0.446667)
 ..controls (4.336575,0.446667) and (4.31792,0.428012)
 ..(4.31792,0.405)\dpicstop
\dpicdraw (4.401253,0.405)
 ..controls (4.401253,0.381988) and (4.419908,0.363333)
 ..(4.44292,0.363333)
 ..controls (4.465932,0.363333) and (4.484587,0.381988)
 ..(4.484587,0.405)\dpicstop
\dpicdraw (4.526253,0.405)
 --(4.651253,0.405)\dpicstop
\draw (4.276253,0.405) node[above left=-2bp]{$ +$};
\draw (4.401253,0.53) node[above=-2bp]{$ v_{\g,\A}$};
\draw (4.526253,0.405) node[above right=-2bp]{$ -$};
\dpicdraw (2.951253,0)
 --(3.151253,0)\dpicstop
\draw (3.151253,0) node[above=-2bp]{$v_{\C,\B}$};
\dpicdraw (3.151253,0)
 --(3.276253,0)\dpicstop
\dpicdraw (3.276253,0)
 --(3.276253,-0.005556)\dpicstop
\dpicdraw (3.276253,0)
 ..controls (3.276253,0.017259) and (3.290244,0.03125)
 ..(3.307503,0.03125)
 ..controls (3.324762,0.03125) and (3.338753,0.017259)
 ..(3.338753,0)\dpicstop
\dpicdraw (3.338753,0)
 --(3.338753,-0.005556)\dpicstop
\dpicdraw (3.338753,0)
 ..controls (3.338753,0.017259) and (3.352744,0.03125)
 ..(3.370003,0.03125)
 ..controls (3.387262,0.03125) and (3.401253,0.017259)
 ..(3.401253,0)\dpicstop
\dpicdraw (3.401253,0)
 --(3.401253,-0.005556)\dpicstop
\dpicdraw (3.401253,0)
 ..controls (3.401253,0.017259) and (3.415244,0.03125)
 ..(3.432503,0.03125)
 ..controls (3.449762,0.03125) and (3.463753,0.017259)
 ..(3.463753,-0)\dpicstop
\dpicdraw (3.463753,-0)
 --(3.463753,-0.005556)\dpicstop
\dpicdraw (3.463753,-0)
 ..controls (3.463753,0.017259) and (3.477744,0.03125)
 ..(3.495003,0.03125)
 ..controls (3.512262,0.03125) and (3.526253,0.017259)
 ..(3.526253,-0)\dpicstop
\dpicdraw (3.526253,-0)
 --(3.526253,-0.005556)\dpicstop
\dpicdraw (3.526253,-0)
 --(3.651253,-0)\dpicstop
\filldraw (3.576253,-0.01875)
 --(3.651253,0)
 --(3.576253,0.01875) --cycle\dpicstop
\dpicdraw (3.628347,0)
 --(3.576253,0)\dpicstop
\draw (3.6023,0) node[above=-2bp]{$ i_b$};
\dpicdraw[fill=white](3.151253,0) circle (0.007874in)\dpicstop
\dpicdraw (3.651253,0)
 --(3.776253,0)
 --(3.797087,0.041667)
 --(3.838753,-0.041667)
 --(3.88042,0.041667)
 --(3.922087,-0.041667)
 --(3.963753,0.041667)
 --(4.00542,-0.041667)
 --(4.026253,0)
 --(4.151253,0)\dpicstop
\dpicdraw (4.151253,0)
 --(4.276253,0)\dpicstop
\dpicdraw (4.401253,0) circle (0.049213in)\dpicstop
\dpicdraw (4.401253,0)
 ..controls (4.401253,0.023012) and (4.382598,0.041667)
 ..(4.359587,0.041667)
 ..controls (4.336575,0.041667) and (4.31792,0.023012)
 ..(4.31792,0)\dpicstop
\dpicdraw (4.401253,0)
 ..controls (4.401253,-0.023012) and (4.419908,-0.041667)
 ..(4.44292,-0.041667)
 ..controls (4.465932,-0.041667) and (4.484587,-0.023012)
 ..(4.484587,0)\dpicstop
\dpicdraw (4.526253,0)
 --(4.651253,0)\dpicstop
\draw (4.276253,0) node[above left=-2bp]{$ +$};
\draw (4.401253,0.125) node[above=-2bp]{$ v_{\g,\B}$};
\draw (4.526253,0) node[above right=-2bp]{$ -$};
\dpicdraw[line width=0.4bp](2.951253,-0.105) circle (0.00109in)\dpicstop
\dpicdraw (2.951253,-0.105)
 --(2.951253,-0.405)\dpicstop
\dpicdraw[line width=0.4bp](2.951253,-0.405) circle (0.00109in)\dpicstop
\dpicdraw (2.951253,-0.405)
 --(3.151253,-0.405)\dpicstop
\dpicdraw[fill=white](3.151253,-0.405) circle (0.007874in)\dpicstop
\draw (3.151253,-0.405) node[above=-2bp]{$v_{\C,\C}$};
\dpicdraw (3.151253,-0.405)
 --(3.276253,-0.405)\dpicstop
\dpicdraw (3.276253,-0.405)
 --(3.276253,-0.410556)\dpicstop
\dpicdraw (3.276253,-0.405)
 ..controls (3.276253,-0.387741) and (3.290244,-0.37375)
 ..(3.307503,-0.37375)
 ..controls (3.324762,-0.37375) and (3.338753,-0.387741)
 ..(3.338753,-0.405)\dpicstop
\dpicdraw (3.338753,-0.405)
 --(3.338753,-0.410556)\dpicstop
\dpicdraw (3.338753,-0.405)
 ..controls (3.338753,-0.387741) and (3.352744,-0.37375)
 ..(3.370003,-0.37375)
 ..controls (3.387262,-0.37375) and (3.401253,-0.387741)
 ..(3.401253,-0.405)\dpicstop
\dpicdraw (3.401253,-0.405)
 --(3.401253,-0.410556)\dpicstop
\dpicdraw (3.401253,-0.405)
 ..controls (3.401253,-0.387741) and (3.415244,-0.37375)
 ..(3.432503,-0.37375)
 ..controls (3.449762,-0.37375) and (3.463753,-0.387741)
 ..(3.463753,-0.405)\dpicstop
\dpicdraw (3.463753,-0.405)
 --(3.463753,-0.410556)\dpicstop
\dpicdraw (3.463753,-0.405)
 ..controls (3.463753,-0.363333) and (3.526253,-0.363333)
 ..(3.526253,-0.405)\dpicstop
\dpicdraw (3.526253,-0.405)
 --(3.526253,-0.410556)\dpicstop
\dpicdraw (3.526253,-0.405)
 --(3.651253,-0.405)\dpicstop
\filldraw (3.576253,-0.42375)
 --(3.651253,-0.405)
 --(3.576253,-0.38625) --cycle\dpicstop
\dpicdraw (3.628347,-0.405)
 --(3.576253,-0.405)\dpicstop
\draw (3.6023,-0.405) node[above=-2bp]{$ i_c$};
\dpicdraw[fill=white](3.151253,-0.405) circle (0.007874in)\dpicstop
\dpicdraw (3.651253,-0.405)
 --(3.776253,-0.405)
 --(3.797087,-0.363333)
 --(3.838753,-0.446667)
 --(3.88042,-0.363333)
 --(3.922087,-0.446667)
 --(3.963753,-0.363333)
 --(4.00542,-0.446667)
 --(4.026253,-0.405)
 --(4.151253,-0.405)\dpicstop
\dpicdraw (4.151253,-0.405)
 --(4.276253,-0.405)\dpicstop
\dpicdraw (4.401253,-0.405) circle (0.049213in)\dpicstop
\dpicdraw (4.401253,-0.405)
 ..controls (4.401253,-0.381988) and (4.382598,-0.363333)
 ..(4.359587,-0.363333)
 ..controls (4.336575,-0.363333) and (4.31792,-0.381988)
 ..(4.31792,-0.405)\dpicstop
\dpicdraw (4.401253,-0.405)
 ..controls (4.401253,-0.428012) and (4.419908,-0.446667)
 ..(4.44292,-0.446667)
 ..controls (4.465932,-0.446667) and (4.484587,-0.428012)
 ..(4.484587,-0.405)\dpicstop
\dpicdraw (4.526253,-0.405)
 --(4.651253,-0.405)\dpicstop
\draw (4.276253,-0.405) node[above left=-2bp]{$ +$};
\draw (4.401253,-0.28) node[above=-2bp]{$ v_{\g,\C}$};
\draw (4.526253,-0.405) node[above right=-2bp]{$ -$};
\dpicdraw (4.651253,0.405)
 --(4.651253,-0.405)\dpicstop
\dpicdraw[fill=black](4.651253,0) circle (0.007874in)\dpicstop
\draw (-0.348747,-0.084) node[left=-2bp]{$\mathrm{N}$};
\end{tikzpicture}%

%% file: Figures/block_diagram.tex
\begin{tikzpicture}[auto, node distance=5.68em, >=stealth', every node/.style={align=center},scale = .3]
        
        % Define Tikz Styles
        \tikzstyle{block} = [draw, thick, fill=blueCol!20, rectangle, inner sep=5pt]
        \tikzstyle{sum} = [draw, circle, node distance=2.84em, inner sep=0pt]
        \tikzstyle{branch} = [circle,inner sep=0pt,minimum size=1mm,fill=black,draw=black]
        \tikzstyle{input} = [coordinate]
        \tikzstyle{output} = [coordinate]
            
            % We start by placing the blocks
            \node [input] (curr) {};
            \node [input, below of = curr, node distance = 2em] (fsw) {};
            \node [block, right of = curr, node distance = 10em, minimum height=12em,yshift = -5em] (mpc) {MPC};
            \node [block, below of=mpc, node distance = 8.5em] (filt) {EST};
            \node [block, below of=filt, node distance = 2em] (delay) {$z^{-1}$};
            \node [branch,right of=mpc, node distance=6em, yshift = 1.42em] (branchU) {};
            \node[coordinate, right of=mpc, node distance=4em, yshift = -1.42em](p){};
            \node[coordinate, left of=mpc, node distance=4em, yshift = -5em](pprev){}; 
            \node[coordinate, left of=mpc, node distance=7em, yshift = -3em](uprev){};
            \node [block,right of=branchU, node distance=5.544em,minimum height=4.26em, minimum width = 4.26em] (motor) {CONVERTER\\SYSTEM};
            \node [block, below of=motor, node distance = 8.52em, xshift=-1.1em] (clk) {$\b{K}$};
            \node [block, below of=motor, node distance = 8.52em, xshift=1.1em] (clk1) {$\b{K}$};
            \node[output, below of=clk, node distance = 1.42em] (xphys1){};
            \node[output, below of=delay, node distance = 3em, left of =xphys1, yshift=-3em] (xphys2){};
            \node[output, left of=mpc, node distance =9em, yshift = -1em] (xphys3){};
            \node[output, below of=clk1, node distance = 1.42em] (Vg1) {};
            \node[output, below of=xphys2, node distance = 3em, left of=Vg1, yshift = -4em] (Vg2) {};
            \node[output, left of=mpc, node distance =10em, yshift = 1em] (Vg3) {};
            
            % Once the nodes are placed, we connect them
            \draw [->] (curr) -- node [pos=0.5] {$\boldsymbol{i}^*(k)$} (curr) -- (curr-|mpc.west);
            \draw [->] (fsw) -- node [pos=0.5] {$f_\sw^*(k)$} (fsw) -- (fsw-|mpc.west);
            \draw [->] (branchU-|mpc.east) -- node [pos=0.55] {$\boldsymbol{u}_\ph(k)$} (branchU) -- (motor);
            \draw [->] (motor.295) -- node {$\b{v}_{\g,\abc}(k)$}(clk1.north-|motor.295);
            \draw [->] (motor.245) -- node [anchor = east]{$\b{i}_\abc(k)$}(clk.north-|motor.245);
            \draw [->] (mpc.east|-p) -- node [pos=0.7] {$\b{p}(k)$} (p) |- (filt.east);
            \draw [->] (filt.west) -|  (pprev) -- node [pos=0.2] {$f_\sw(k)$}(pprev-|mpc.west);
            \draw [->] (branchU) |- (delay);
            \draw [->] (delay) -| (uprev) -- node[pos = 0.5] {$\boldsymbol{u}_\ph(k-1)$}(uprev-|mpc.west);
            \draw [->] (clk.south) -- (xphys1) |- node [pos =0.1,anchor=east]{$\boldsymbol{x}_\ph(k)$}(xphys2) -| (xphys3) -- node [pos = 0.75]{$\boldsymbol{x}_\ph(k)$}(xphys3-|mpc.west);
            \draw [->] (clk1.south) -- (Vg1) |- node [pos =0.1]{$\boldsymbol{v}_\g(k)$}(Vg2) -| (Vg3) -- node [pos = 0.81]{$\boldsymbol{v}_\g(k)$}(Vg3-|mpc.west);
            \draw[thick, dotted] ($(mpc.north) + (0.0,0.3)$) -| ($(branchU) + (0.3,0.0)$) |- ($(delay.south) + (0.0, -1.5)$) -| node[pos = 0.0, anchor = south east]{CONTROLLER}($(uprev.west) + (-.3,0.0)$) |- ($(mpc.north) + (0.0,0.3)$);
            
\end{tikzpicture}

%% file: Figures/f_sw-ss.tex
% Use this code to include graph into document:
%
%\begin{figure}
%	\setlength{\fwidth}{0.6\linewidth}
%	\setlength{\fheight}{0.8\fwidth}
%	\input{Images/opp-cur-d5.tex}
%\end{figure}
{\scriptsize
\begin{tikzpicture}
\def\dx{-0.7};
\def\dy{-0.4};
\begin{axis}[
grid = major,
xmin=.5,
xmax=1.5,
xlabel = {Time [\si{\second}]},
ylabel = {$f_\sw$ [\si{\hertz}]},
ylabel near ticks,
width=\fwidth,
height=\fheight,
legend pos=north east,
ymax = 325,
ymin = 205,
]
\addplot [restrict x to domain=0:3, color=black, dashed] table[x=t, y=f_sw_ref,mark=none] {Figures/f_sw-ss.txt};
\addplot [restrict x to domain=0:3,color=green, thick] table[x=t, y=f_sw_2,mark=none] {Figures/f_sw-ss.txt};
\addplot [restrict x to domain=0:3,color=blue] table[x=t, y=f_sw_3,mark=none] {Figures/f_sw-ss.txt};
\legend{ref, FT, FL}
\end{axis}
\end{tikzpicture}}

%% file: Figures/i-ss.tex
% Use this code to include graph into document:
%
%\begin{figure}
%	\setlength{\fwidth}{0.6\linewidth}
%	\setlength{\fheight}{0.8\fwidth}
%	\input{Figures/opp-cur-d5.tex}
%\end{figure}
{\scriptsize
\begin{tikzpicture}
\def\dx{-0.7};
\def\dy{-0.4};
\begin{axis}[
grid = major,
xticklabel=\pgfkeys{/pgf/number format/.cd,fixed,precision=3,zerofill}\pgfmathprintnumber{\tick},
xmin=1,
xmax=1.02,
xlabel = {Time [\si{\second}]},
ylabel = {$\b{i}$ [\si{\pu}]},
ylabel near ticks,
width=\fwidth,
height=\fheight,
% legend pos=south west,
legend style={at={(.8,0.025)},anchor=south west},
ymax = 1.2,
ymin = -1.2,
]

\addplot [restrict x to domain=1:1.02, color=black, dashed] table[x=t, y=i_ref_1, mark=none] {Figures/i-ss.txt};
\addplot [restrict x to domain=1:1.02, color=green] table[x=t, y=i_ctrl2_1, mark=none] {Figures/i-ss.txt};
\addplot [restrict x to domain=1:1.02, color=blue] table[x=t, y=i_ctrl3_1, mark=none] {Figures/i-ss.txt};

\addplot [restrict x to domain=1:1.02, color=black, dashed] table[x=t, y=i_ref_2, mark=none] {Figures/i-ss.txt};
\addplot [restrict x to domain=1:1.02, color=green] table[x=t, y=i_ctrl2_2, mark=none] {Figures/i-ss.txt};
\addplot [restrict x to domain=1:1.02, color=blue] table[x=t, y=i_ctrl3_2, mark=none] {Figures/i-ss.txt};

\addplot [restrict x to domain=1:1.02, color=black, dashed] table[x=t, y=i_ref_3, mark=none] {Figures/i-ss.txt};
\addplot [restrict x to domain=1:1.02, color=green] table[x=t, y=i_ctrl2_3, mark=none] {Figures/i-ss.txt};
\addplot [restrict x to domain=1:1.02, color=blue] table[x=t, y=i_ctrl3_3, mark=none] {Figures/i-ss.txt};

\legend{ref, FT, FL}
\end{axis}
\end{tikzpicture}}

%% file: Figures/i-iramp.tex
% Use this code to include graph into document:
%
%\begin{figure}
%	\setlength{\fwidth}{0.6\linewidth}
%	\setlength{\fheight}{0.8\fwidth}
%	\input{Figures/opp-cur-d5.tex}
%\end{figure}
{\scriptsize
\begin{tikzpicture}
\def\dx{-0.7};
\def\dy{-0.4};
\begin{axis}[
grid = major,
xticklabel=\pgfkeys{/pgf/number format/.cd,fixed,precision=3,zerofill}\pgfmathprintnumber{\tick},
xmin=1.205,
xmax=1.225,
xlabel = {Time [\si{\second}]},
ylabel = {$\b{i}$ [\si{\pu}]},
ylabel near ticks,
width=\fwidth,
height=\fheight,
% legend pos=south west,
legend style={at={(.8,0.025)},anchor=south west},
ymax = 1.2,
ymin = -1.2,
]

\addplot [restrict x to domain=1.205:1.225, color=black, dashed] table[x=t, y=i_ref_1, mark=none] {Figures/ramp.txt};
\addplot [restrict x to domain=1.205:1.225, color=green] table[x=t, y=i_ctrl2_1, mark=none] {Figures/ramp.txt};
\addplot [restrict x to domain=1.205:1.225, color=blue] table[x=t, y=i_ctrl3_1, mark=none] {Figures/ramp.txt};

\addplot [restrict x to domain=1.205:1.225, color=black, dashed] table[x=t, y=i_ref_2, mark=none] {Figures/ramp.txt};
\addplot [restrict x to domain=1.205:1.225, color=green] table[x=t, y=i_ctrl2_2, mark=none] {Figures/ramp.txt};
\addplot [restrict x to domain=1.205:1.225, color=blue] table[x=t, y=i_ctrl3_2, mark=none] {Figures/ramp.txt};

\addplot [restrict x to domain=1.205:1.225, color=black, dashed] table[x=t, y=i_ref_3, mark=none] {Figures/ramp.txt};
\addplot [restrict x to domain=1.205:1.225, color=green] table[x=t, y=i_ctrl2_3, mark=none] {Figures/ramp.txt};
\addplot [restrict x to domain=1.205:1.225, color=blue] table[x=t, y=i_ctrl3_3, mark=none] {Figures/ramp.txt};

\legend{ref, FT, FL}
\end{axis}
\end{tikzpicture}}

%% file: Figures/i-istep.tex
% Use this code to include graph into document:
%
%\begin{figure}
%	\setlength{\fwidth}{0.6\linewidth}
%	\setlength{\fheight}{0.8\fwidth}
%	\input{Figures/opp-cur-d5.tex}
%\end{figure}
{\scriptsize
\begin{tikzpicture}
\def\dx{-0.7};
\def\dy{-0.4};
\begin{axis}[
grid = major,
xticklabel=\pgfkeys{/pgf/number format/.cd,fixed,precision=3,zerofill}\pgfmathprintnumber{\tick},
xmin=.595,
xmax=.615,
xlabel = {Time [\si{\second}]},
ylabel = {$\b{i}$ [\si{\pu}]},
ylabel near ticks,
width=\fwidth,
height=\fheight,
legend pos=south west,
ymax = 1.2,
ymin = -1.2,
]
\addplot [restrict x to domain=.58:.62, color=black, dashed] table[x=t, y=i_ref_1, mark=none] {Figures/step.txt};
\addplot [restrict x to domain=.58:.62, color=green] table[x=t, y=i_ctrl2_1, mark=none] {Figures/step.txt};
\addplot [restrict x to domain=.58:.62, color=blue] table[x=t, y=i_ctrl3_1, mark=none] {Figures/step.txt};

\addplot [restrict x to domain=.58:.62, color=black, dashed] table[x=t, y=i_ref_2, mark=none] {Figures/step.txt};
\addplot [restrict x to domain=.58:.62, color=green] table[x=t, y=i_ctrl2_2, mark=none] {Figures/step.txt};
\addplot [restrict x to domain=.58:.62, color=blue] table[x=t, y=i_ctrl3_2, mark=none] {Figures/step.txt};

\addplot [restrict x to domain=.58:.62, color=black, dashed] table[x=t, y=i_ref_3, mark=none] {Figures/step.txt};
\addplot [restrict x to domain=.58:.62, color=green] table[x=t, y=i_ctrl2_3, mark=none] {Figures/step.txt};
\addplot [restrict x to domain=.58:.62, color=blue] table[x=t, y=i_ctrl3_3, mark=none] {Figures/step.txt};

\legend{ref, FT, FL}
\end{axis}
\end{tikzpicture}}

%% file: Figures/fsw-fswstep.tex
% Use this code to include graph into document:
%
%\begin{figure}
%	\setlength{\fwidth}{0.6\linewidth}
%	\setlength{\fheight}{0.8\fwidth}
%	\input{Images/opp-cur-d5.tex}
%\end{figure}
{\scriptsize
\begin{tikzpicture}
\def\dx{-0.7};
\def\dy{-0.4};
\begin{axis}[
grid = major,
xmin=0.35,
xmax=0.5,
xlabel = {Time [\si{\second}]},
ylabel = {$f_\sw$ [\si{\hertz}]},
ylabel near ticks,
width=\fwidth,
height=\fheight,
legend pos=south east,
ymax = 325,
ymin = 205,
]
\addplot [restrict x to domain=0:3, color=black, dashed] table[x=t, y=f_sw_ref,mark=none] {Figures/f_sw-step.txt};
\addplot [restrict x to domain=0:3,color=green, thick] table[x=t, y=f_sw_2,mark=none] {Figures/f_sw-step.txt};
\addplot [restrict x to domain=0:3,color=blue] table[x=t, y=f_sw_3,mark=none] {Figures/f_sw-step.txt};
\legend{ref, FT, FL}
\end{axis}
\end{tikzpicture}}

%% file: Figures/i-fswstep.tex
% Use this code to include graph into document:
%
%\begin{figure}
%	\setlength{\fwidth}{0.6\linewidth}
%	\setlength{\fheight}{0.8\fwidth}
%	\input{Figures/opp-cur-d5.tex}
%\end{figure}
{\scriptsize
\begin{tikzpicture}
\def\dx{-0.7};
\def\dy{-0.4};
\begin{axis}[
grid = major,
xticklabel=\pgfkeys{/pgf/number format/.cd,fixed,precision=3,zerofill}\pgfmathprintnumber{\tick},
xmin=.39,
xmax=.41,
xlabel = {Time [\si{\second}]},
ylabel = {$\b{i}$ [\si{\pu}]},
ylabel near ticks,
width=\fwidth,
height=\fheight,
legend pos=south west,
ymax = 1.2,
ymin = -1.2,
]
\addplot [restrict x to domain=.39:.41, color=black, dashed] table[x=t, y=i_ref_1, mark=none] {Figures/i-fswstep.txt};
\addplot [restrict x to domain=.39:.41, color=green] table[x=t, y=i_ctrl2_1, mark=none] {Figures/i-fswstep.txt};
\addplot [restrict x to domain=.39:.41, color=blue] table[x=t, y=i_ctrl3_1, mark=none] {Figures/i-fswstep.txt};

\addplot [restrict x to domain=.39:.41, color=black, dashed] table[x=t, y=i_ref_2, mark=none] {Figures/i-fswstep.txt};
\addplot [restrict x to domain=.39:.41, color=green] table[x=t, y=i_ctrl2_2, mark=none] {Figures/i-fswstep.txt};
\addplot [restrict x to domain=.39:.41, color=blue] table[x=t, y=i_ctrl3_2, mark=none] {Figures/i-fswstep.txt};

\addplot [restrict x to domain=.39:.41, color=black, dashed] table[x=t, y=i_ref_3, mark=none] {Figures/i-fswstep.txt};
\addplot [restrict x to domain=.39:.41, color=green] table[x=t, y=i_ctrl2_3, mark=none] {Figures/i-fswstep.txt};
\addplot [restrict x to domain=.39:.41, color=blue] table[x=t, y=i_ctrl3_3, mark=none] {Figures/i-fswstep.txt};

\legend{ref, FT, FL}
\end{axis}
\end{tikzpicture}}